\documentclass[11pt]{article}
\usepackage{epsfig}
\usepackage{amssymb}
\usepackage{amsmath}

\newtheorem{theorem}{Theorem}[section]
\newtheorem{lemma}[theorem]{Lemma}
\newtheorem{proposition}[theorem]{Proposition}

\newtheorem{corollary}[theorem]{Corollary}
\newtheorem{remark}[theorem]{Remark}

\newtheorem{problem}[theorem]{Problem}
\newtheorem{definition}[theorem]{Definition}
\newtheorem{example}[theorem]{Example}

\newtheorem{conjecture}[theorem]{Conjecture}

\newcommand{\C}{\mathbf{C}}

\newenvironment{proof}{{\it Proof:\/}}{$\Box$\vskip 0.08in}
\newcommand\et{
\genfrac{}{}{0pt}{}
{
\mbox{
$\genfrac{}{}{0pt}{}
{}{
\genfrac{}{}{0pt}{}
{\doteq} {I_t}
}
$}
}
{}
}

\newcommand\ea{
\genfrac{}{}{0pt}{}
{
\mbox{
$\genfrac{}{}{0pt}{}
{}{
\genfrac{}{}{0pt}{}
{\doteq} {I_A}
}
$}
}
{}
}

\newcommand\eax{
\genfrac{}{}{0pt}{}
{
\mbox{
$\genfrac{}{}{0pt}{}
{}{
\genfrac{}{}{0pt}{}
{\doteq} {I_{a,x}}
}
$}
}
{}
}

\newcommand\eaa{
\genfrac{}{}{0pt}{}
{
\mbox{
$\genfrac{}{}{0pt}{}
{}{
\genfrac{}{}{0pt}{}
{\doteq} {I_{a,a+a^{-1}}}
}
$}
}
{}
}

\newcommand\sg{
\genfrac{}{}{0pt}{}{\smile}{\frown}
}
\begin{document}
\begin{center}
\bigskip

{\LARGE \baselineskip=10pt {\ 5-move equivalence classes of links and
their algebraic invariants}}
\end{center}
\centerline{\bf Paper dedicated to Lou Kauffman on his 60th birthday}
\begin{tabular}{lll}
Mieczys{\l}aw K. D{\c a}bkowski & Makiko Ishiwata
 &J\'ozef H. Przytycki
\\
Department of Mathematics & Osaka City University
& Department of Mathematics \\
University of Texas at Dallas & OCAMI
&George Washington University \\
Richardson, TX 75083 & Osaka-city 558-8585
& On sabbatical leave at UMD \\
 mdab@utdallas.edu &  misiwata@gmail.com & przytyck@gwu.edu\\
\end{tabular}
\ \\
Key Words: {\it 5-moves, rational moves, links, knots, Jones Polynomial,
Kauffman polynomial, Burnside group, Pretzel link, Montesinos link}\\
{Mathematics Subject Classification 2000: 57M25, 57M27}\ \\

\vspace{27pt} Abstract.\newline
{\footnotesize {\
We start a systematic analysis of links up to 5-move equivalence.
Our motivation is to develop tools which later can be used to
study skein modules based on the skein relation being deformation of a
5-move (in an analogous way as the Kauffman skein module is a
deformation of a 2-move, i.e. a crossing change). Our main tools are
Jones and Kauffman polynomials and the fundamental group of the
2-fold branch cover of $S^3$ along a link. We use also the fact that
a 5-move is a composition of two rational $\pm (2,2)$-moves (i.e.
$\pm \frac{5}{2}$-moves) and rational moves can be analyzed using
the group of Fox colorings and its non-abelian version,
the Burnside group of a link.
One curious observation is that links related by one $(2,2)$-move are not
$5$-move equivalent. In particular, we partially classify (up to 5-moves)
 3-braids, pretzel and Montesinos links, and links up to 9 crossings.}}

\ \newline
\tableofcontents

\section{Introduction\label{1}}

A tangle move is a local modification of a link in which a tangle
$T_{A}$ is replaced by a tangle $T_{B}$, Fig. 1.1. \\
\ \\
\centerline{\psfig{figure=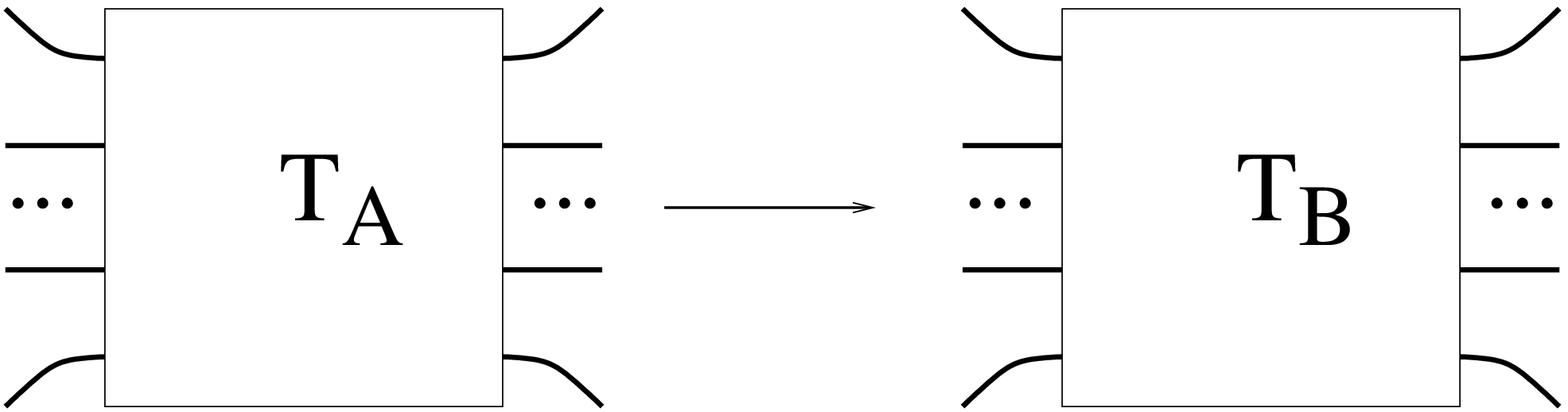,height=2.1cm}} 
\begin{center}
Fig. 1.1; Tangle move
\end{center}

Our interest in
tangle moves on links has been motivated by our analysis of skein modules of
3-dimensional manifolds. Skein relations for links might be viewed as
deformations of tangle moves. The simplest moves that reduce every
link in $S^{3}$ into a trivial link are a smoothing of a crossing and
a crossing change. A deformation of a smoothing leads to
Kauffman bracket skein module and a deformation of a crossing change
leads to, in the oriented case,  Jones and Homflypt
skein modules, and in the unoriented case, Kauffman skein module.
In the last case the deformation is of
the form $L_{+} + L_{-}=xL_{0}+ xL_{\infty }$ ( Fig. 1.2).
\ \newline
\ \newline
\centerline{\psfig{figure=L+L-L0Linf.eps,height=1.9cm}} 

\begin{center}
Fig. 1.2
\end{center}

If a move is an unlinking move (i.e. every link can be reduced to a trivial
link) then some deformations of the move can lead to a skein module of $S^3$
generated by trivial links. This is the case for the Kauffman bracket,
Homflypt and Kauffman skein modules (see \cite{H-P} or \cite{Pr-2,Pr-7,Pr-8}
for a survey of skein modules).
A 3-move (\parbox{2.7cm}{\ \psfig{figure=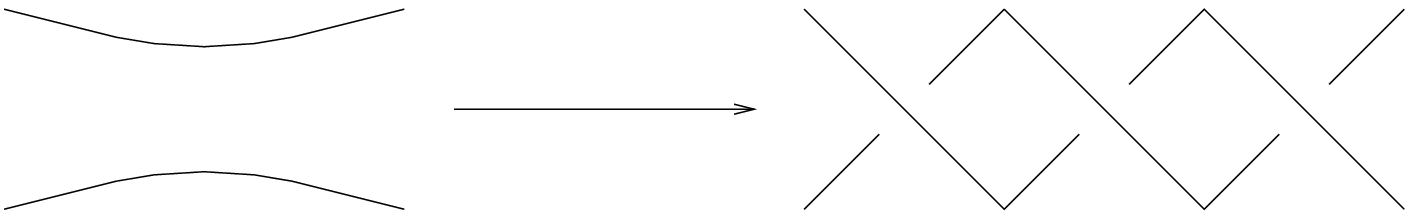,height=0.4cm}} ) is
probably the simplest move after the crossing change.
For over 20 years it was an open problem (the Montesinos-Nakanishi conjecture)
as to whether or not every link can be reduced to a trivial link via 3-moves.
We finally disproved it in 2002 \cite{D-P-1}.
A 4-move (\parbox{3.6cm}{\psfig{figure=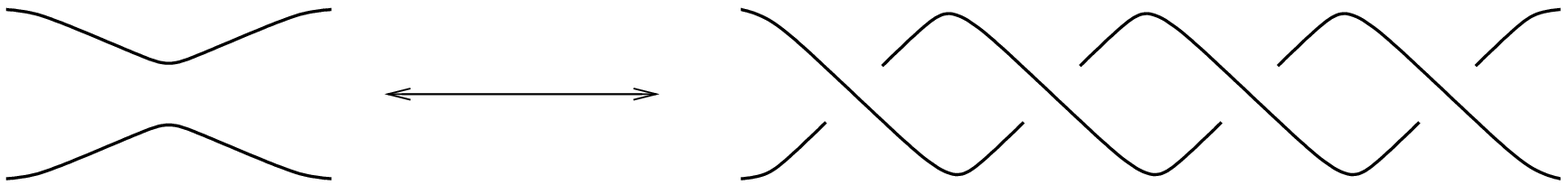,height=0.4cm}})
preserves the number of components  of a link so it makes sense to
study $4$-moves on knots separately. The Nakanishi conjecture,
formulated in 1979, stated
that every knot can be unknotted via 4-moves.
This conjecture remains still open. However, the related question (of Kawauchi)
 for links of three or more components\footnote{For links of two
components the Kawauchi question has the form: can any 2-component
link be reduced by $4$-moves to the trivial link of two components, $T_2$,
or the Hopf link, $H$? The problem is not solved yet.} has been
settled in \cite{D-P-2}.
It is easy to show that not every link is $5$-move equivalent
to a trivial link. For example, the Jones polynomial can be used to demonstrate
that the figure eight knot ($4_1$ in \cite{Rol}) is not
5-move equivalent to any trivial link \cite{Pr-1}. We will develop methods
of analyzing 5-moves using the Jones and Kauffman polynomials in Sections 3
and 4 (compare \cite{Pr-1}).
One can introduce a more delicate move, called $(2,2)$-move
(\parbox{1.6cm}{\psfig{figure=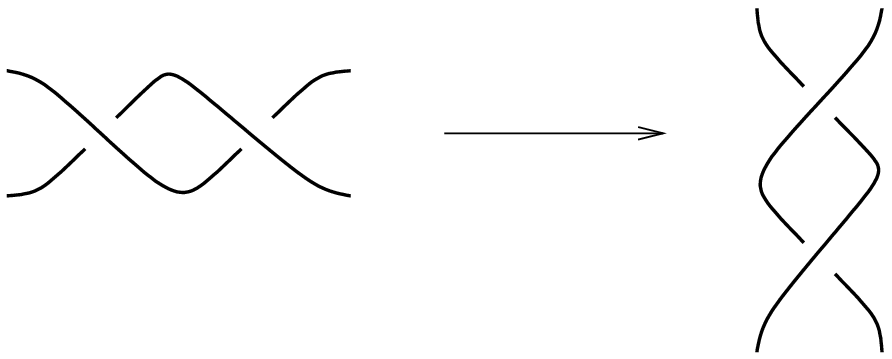,height=0.6cm}})
such that a $5$- move is a combination of a $(2,2)$-move and its mirror
image $(-2,-2)$-move
(\parbox{1.8cm}{\psfig{figure=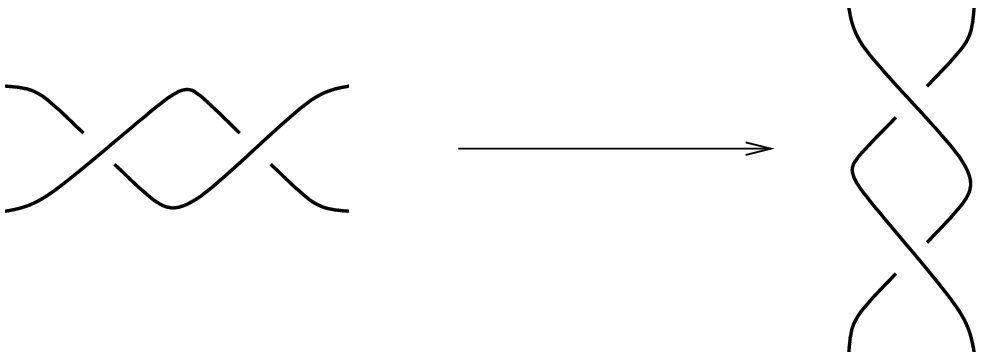,height=0.6cm}}), as
illustrated in Figure 1.3 \cite{H-U,Pr-3}. \newline

\centerline{\psfig{figure=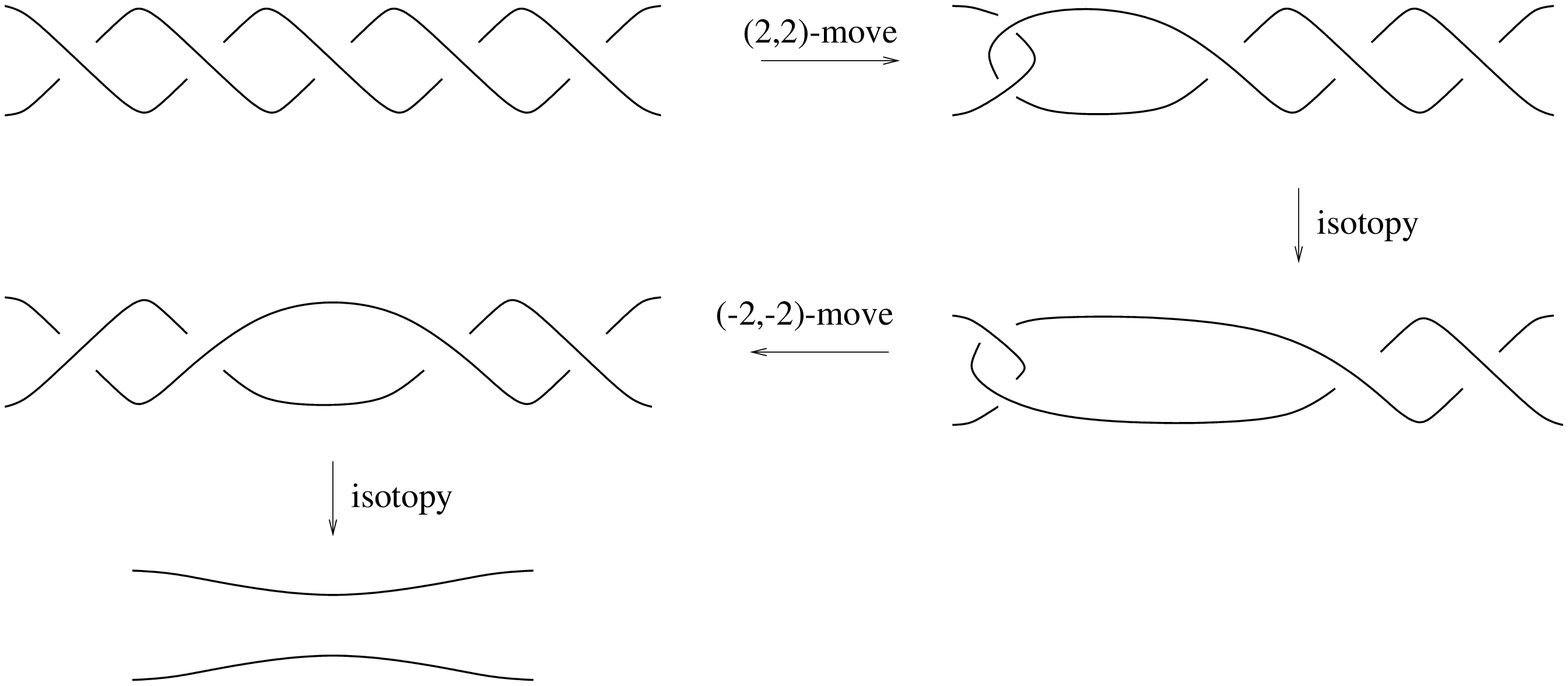,height=4.9cm}}
\centerline{Figure 1.3}

The Harikae-Nakanishi-Uchida conjecture, formulated in 1992,
states that every link can be reduced to a trivial link via $\pm(2,2)$-moves.
This conjecture was disproved in \cite{D-P-2}.
One can try to find $(2,2)$-move equivalence classes of
links. The main objects of this paper are links up to 5-moves, but
because a 5-move is
a combination of $\pm(2,2)$-moves we devote the first two sections of
the paper to the analysis of links up to $\pm(2,2)$-moves, in particular,
algebraic links, 3-braid links, and links up to $9$ crossings. \\

The paper is organized as follows: we introduce gradually
invariants of $(2,2)$- and $5$-moves and we illustrate
constructed invariants analyzing
some family of links (e.g. rational links or algebraic links).
Finally we use all our invariants to (partially)
 classify $5$-move equivalences of $3$-braids, Montesinos links,
 and links up to $9$ crossings.
\\
\ \\
\section{Invariants of $(2,2)$-moves and their applications}
~~\par
We discuss in this section invariants of links which are
preserved by $(2,2)$- or 5-moves.
The simplest of such invariants is the space of Fox $5$-colorings,
${Col}_5(L)$. We describe its use in the next subsection.

\subsection{Fox $n$-colorings and algebraic tangles}\label{Subsection 2.1}
~~\par
The first invariant we apply to analyze rational moves is the group
of Fox $n$-colorings.
We recall first the notion of a rational $\frac{n}{m}$-move and
$n$-rational-equivalence of links and tangles.

\begin{definition}\label{Definition 2.1}
(i) Rational $\frac{n}{m}$-move is a tangle move in which
the $[0]$-tangle is replaced by $[\frac{n}{m}]$-tangle (see Fig.2.1
for rational $\frac{5}{2}$- and $-\frac{5}{2}$-moves.)\\
(ii) We say that two links (or tangles) $L$ and $L'$
are $n$-rationally-equivalent if $L'$ can be
obtained from $L$ by a finite number
of rational $\frac{ns}{m}$-moves ($m$ and $s$ are any non-zero integers).
\end{definition}

\centerline{\psfig{figure=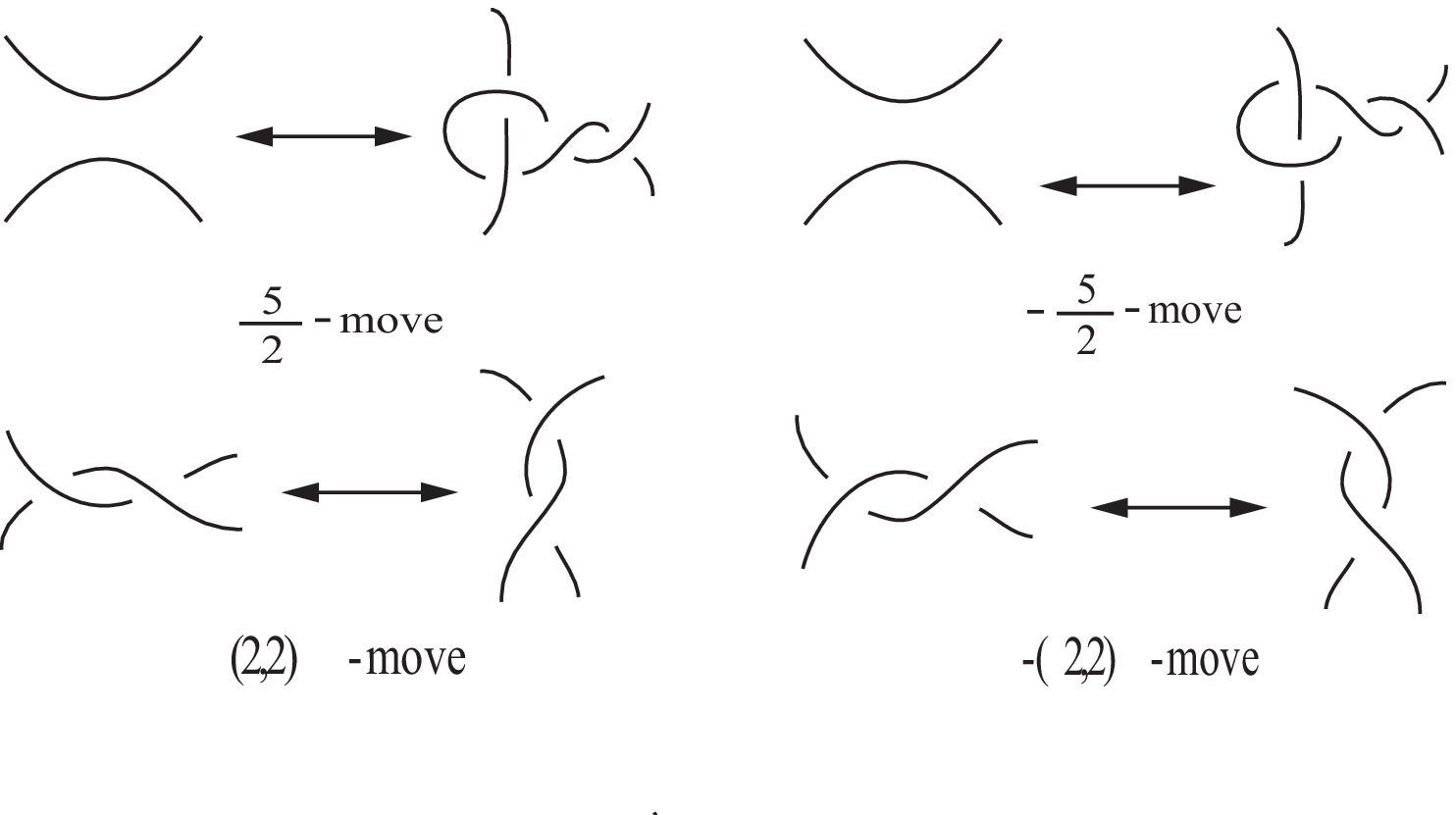,height=3.6cm}}
\centerline{Figure 2.1}

We noted in \cite{Pr-6} that $5$-rationally-equivalence is
the same as $\frac{5}{2}$-move equivalence which in turn is
the same as $(2,2)$-move equivalence in which
we allow the finite number of $\pm (2,2)$-moves (compare Figure 2.1).

Recall that the group of Fox $n$-colorings of a link $L$, ${Col}_n(L)$,
satisfies
${Col}_n(L)=H_1(M_{L}^{(2)}; Z_n) \oplus Z_n$,
where $M_L^{(2)}$ denotes the double branched cover over $S^3$ along $L$
(see \cite{Pr-3} for the combinatorial definition and detailed discussion).

\begin{lemma}\label{Lemma 2.2}
${Col}_n(L)$
is preserved by a rational $\frac{ns}{m}$-move for any non-zero $m$ and $s$.
In particular, $Col_n(L)$ is preserved by $n$-moves.
\end{lemma}

For a trivial link of $k$ components, $T_k$, we have ${Col}_n(T_k)=Z_n^k$.

${Col}_n(L)$ is a rather weak invariant of links but it can be used
as the first step in classifying links up
to $\frac{ns}{m}$-moves ($n$-rational-equivalence).

If $n$ is a prime number then ${Col}_n(L)$ brings the same information
as its order, which we denote by ${col}_n(L)$.


We will give a few applications of Fox $n$-colorings.
We use standard Conway notation for rational tangles\footnote{Our notation
follows Conway's \cite{Con} and agrees with that of Kawauchi book
\cite{Kaw}, but the mirror image notation is often in use \cite{K-L}.}
 (compare Fig. 2.4)
and for the numerator $T^N$
(\ \parbox{1.1cm}{\psfig{figure=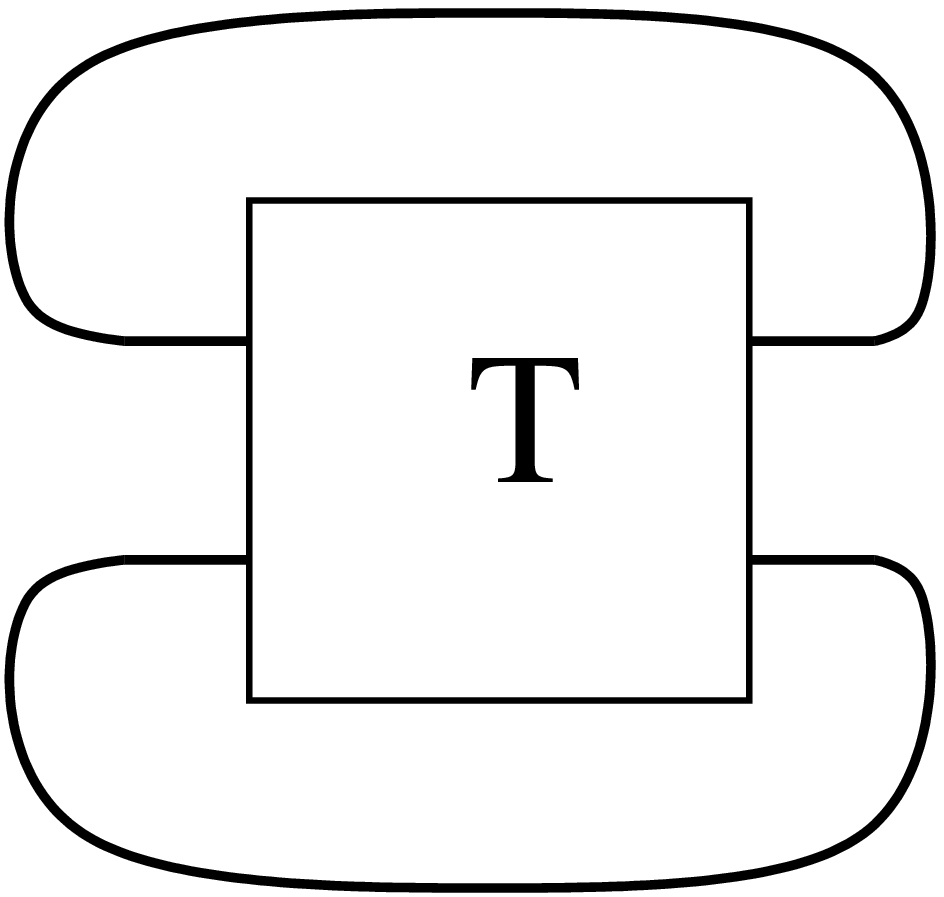,height=0.8cm}}),
and for the denominator $T^D$
(\ \parbox{1.6cm}{\psfig{figure=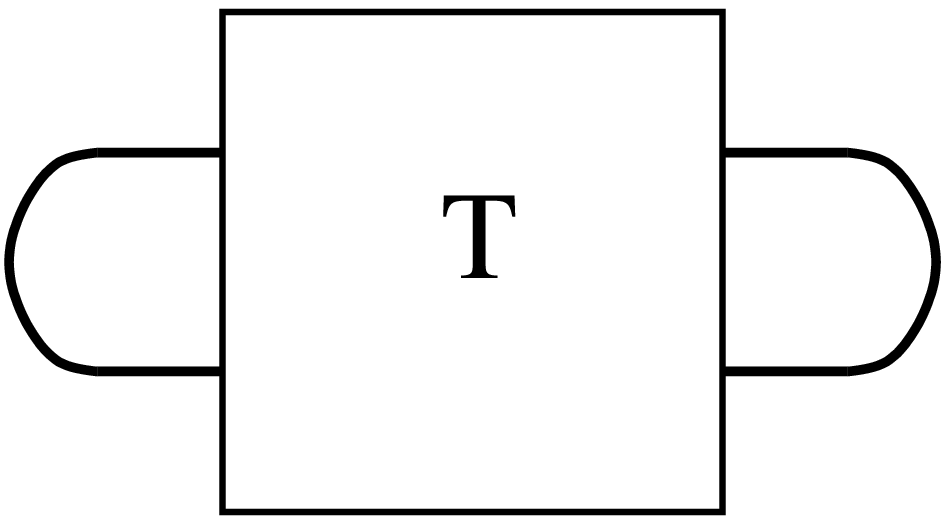,height=0.8cm}})
of a tangle $T$, and for the product of two tangles $T_A*T_B$
(\ \parbox{2.5cm}{\psfig{figure=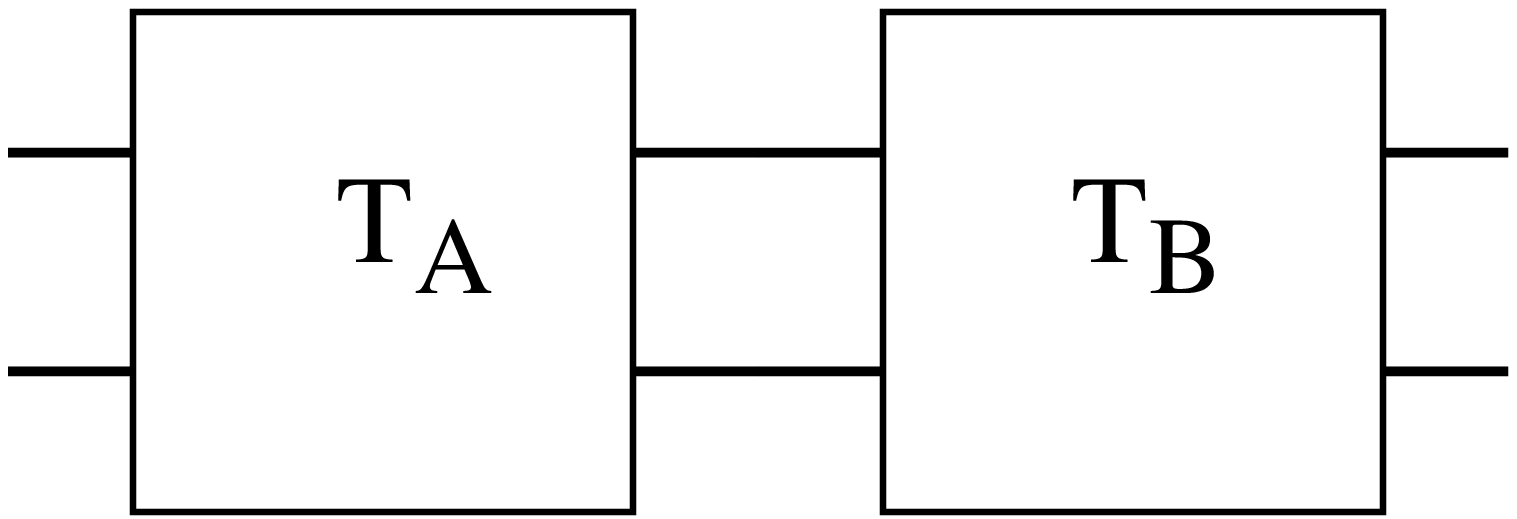,height=0.8cm}}).

\begin{lemma}\label{Lemma 2.3}
The $2$-tangles $[\infty]$, $[0],[1],...,[n-1]$ are pairwise not
$n$-rationally-equivalent, in particular the $2$-tangles
$[\infty]$, $[0],[-1],[1],[-2],$ and $[2]$ are pairwise not
$(2,2)$-move equivalent.
\end{lemma}

\begin{proof}
One easily checks that ${col}_n([k]^N)=n^2$ if and only if
$k$ is a multiple of $n$.
Therefore for any $i$, $0 \leq i \leq n-1$, or $i=\infty$ there exist
exactly one $j$, $(0 \leq j \leq n-1$, or $j=\infty$),
such that ${col}_n(([i] \ast [j])^N)=n^2$: of course
for $i=\infty$ one has $j=\infty$ and for $0 \leq i \leq n-1$ one has $j=n-i$.
From this follows that no pair of elements from $[\infty], [0],[1],...,[n-1]$
are $n$-rationally-equivalent because if $[i]$ and $[i']$ would
be $n$-rationally equivalent then for any
$j$, ${col}_n(([i] \ast [j])^N)={col}_n(([i'] \ast [j])^N)$
which contradicts the previous conclusion\footnote{One can formulate
Lemma 2.3 in more sophisticated language: all tangles
$[\infty]$, $[0],[1],...,[n-1]$ represent different Lagrangians
in a symplectic space $Z_n^2$, see \cite{DJP,Pr-6}.}.
\end{proof}

We proved in \cite{Pr-6}  that any algebraic tangle\footnote{Algebraic
tangles were introduced by Conway in \cite{Con}. They are obtained from
2-tangles of no more than one crossing, by product and rotation
operations. They have a natural generalization to $n$-tangles, in which
case they are called $n$-algebraic tangles \cite{P-Ts}.}
 is $n$-rationally-equivalent
 to one of $n+1$ tangles of Lemma 2.3.
In the case of $n=5$, reduction is a pleasure exercise, see
 \cite{DIP}.
In Subsection 2.5, we demonstrate similar results for
$5$-moves and rational tangles.
Let us now  put $n=5$ and illustrate Lemma 2.2 by another example used
later in classification of Montesinos
links up to $5$-moves.

\begin{example}\label{Example 2.4}
Consider links $L(T_A,k) =
(T_A*\underbrace{[\frac{2}{5}]*...*[\frac{2}{5}]}_{k\ times})^N$
as illustrated in Figure 2.2. Then for $k \geq 1$ we have
$col_5(L(T_A,k))=5^{k-1} col_5(T_A^D)$, and therefore these links
 represent pairwise different $(2,2)$-move equivalence classes.
To see this notice that the rational $\frac{2}{5}$ tangle can be changed
by a $(2,2)$-move to $\infty$ tangle $()()$.
Therefore $L(T_A,k)$ is (2,2)-move equivalent
to $T_A^D \sqcup \underbrace{T_1   \sqcup ...\sqcup T_1}_{(k-1)\ times}\
(k \geq 1),$ (Fig.2.2),
for which we easily count the number of $5$-colorings.
\end{example}

\centerline{\psfig{figure=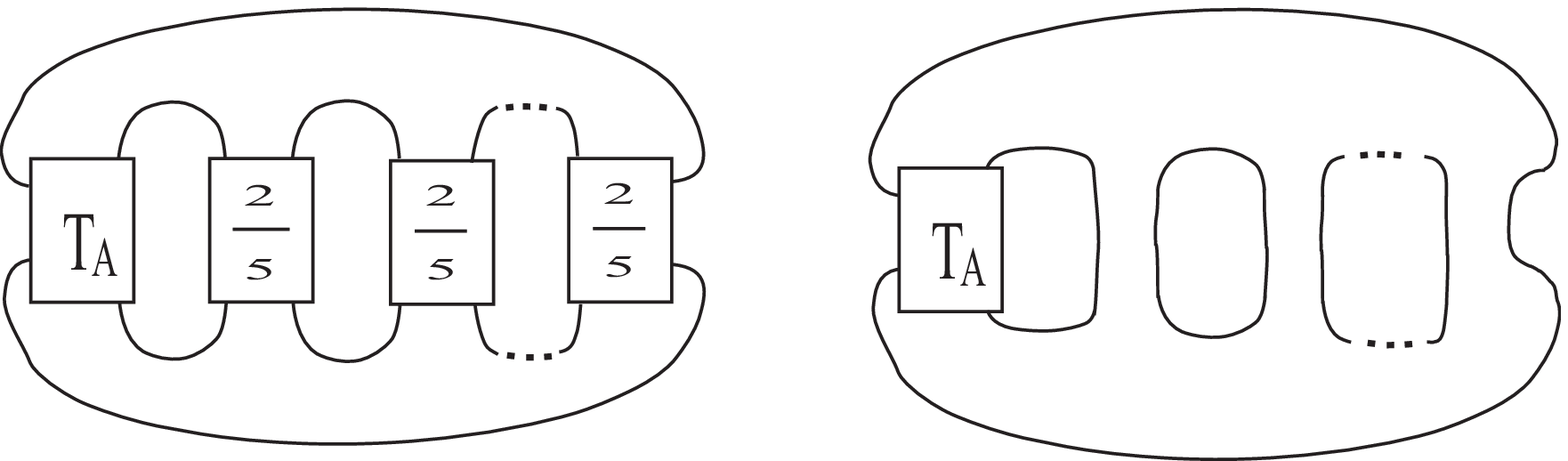,height=3cm}}
\centerline{Figure 2.2}

Notice that $col_5(L)$ cannot distinguish $L$ from the trivial link of
$log_5(col_5(L)$) components.

\subsection{Burnside group of links}

The group of Fox $n$-colorings can be generalized to its non-abelian
version, the $n$th Burnside group of links, $B_n(L)$. This group, which is
also preserved by $\frac{ns}{m}$-rational moves on links, was introduced
in \cite{D-P-1} and used to disprove the Harikae-Nakanishi-Uchida conjecture,
in particular to show that the knots $9_{40}$ and $9_{49}$ are not
$(2,2)$-move equivalent to trivial links.
Recall that the $n$th Burnside group of links, satisfies
$B_n(L)={\pi}_1(M_L^{(2)})/(w^n)$  where the subgroup $(w^n)$ is normally
generated by all elements $w^n$,  $w\in {\pi}_1(M_L^{(2)})$.

\subsection{$(2,2)$-move equivalence classes of algebraic links, 3-braids,
and links up to 9 crossings}\label{Section 2.3}
In this subsection, we summarize and slightly improve the result in
\cite{DIP,Pr-5,Pr-6} (we observe that the knot $9_{49}$ is
related by one 5-move to the mirror image $\bar 9_{49}$).

\begin{theorem}\label{Theorem 2.5}
\begin{enumerate}
\item[(i)] The knots $9_{40}$ and $9_{49}$ are not $(2,2)$-move
equivalent to trivial links. Thus the Harikae-Nakanishi-Uchida conjecture
does not hold.
\item[(ii)] Every algebraic link (in the Conway sense)
 is $(2,2)$-move equivalent to a trivial link.
\item[(iii)]
Every link up to $9$ crossings\footnote{We were informed by S.Jablan that
he checked that every prime link up to 11 crossings and every prime knot up
to 12 crossings is $(2,2)$-move equivalent to a
trivial link or to one of the knots $9_{40}$,
or $9_{49}$ (or their mirror images). Among prime alternating links of
12 crossings there are 3 undecided cases, the links $12_{3*},12_{4*}$, and
$12_{7*}$ \cite{J-S} (in Caudron list of basic polyhedra the names
12C, 12D, and 12G are used and in the Jablan-Sazdanovic book, the corresponding
links are illustrated in Fig.1.74).} is $(2,2)$-move equivalent to a
trivial link or to one of the knots $9_{40}$, its mirror image $\bar 9_{40}$,
or $9_{49}$.
\item[(iv)]
Every closed $3$-braid is $(2,2)$-move equivalent to a trivial
link or to the closure of the braids $(\sigma_1\sigma_2)^6$,
$(\sigma_1\sigma_2)^{12}$ or $(\sigma_1\sigma_2)^{-12}$.

\end{enumerate}
\end{theorem}
\begin{proof}
Part (i) has been proven in \cite{D-P-2} using the fifth Burnside groups
of links.\\
Part (ii) has been demonstrated in \cite{DIP} (compare Lemma 2.10).\\
(iii) It has been demonstrated in \cite{DIP} that any link
up to 9 crossing is $(2,2)$-move equivalent to $9_{40}$, $9_{49}$ or
their mirror images. The proof uses case by case analysis
of non-algebraic links which have at most 9 crossings. The list, which
we will use later, is as follows (up to mirror image):
$8_{18}$, $9_{34}$, $9_{39}$, $9_{40}$, $9_{41}$,
$9_{47}$, $9_{49}$, $9_{40}^{2}$,$9_{41}^{2}$,$9_{42}^{2}$,$9_{61}^{2}$.
The Burnside group argument shows that the
links $9_{40}$, $9_{49}$, $9^2_{40}$, $9_{61}^{2}$ are not $(2,2)$-move
equivalent to trivial links.
We also noticed, \cite{D-P-2,DIP,Pr-6}, that $9^2_{40}$ and $\bar 9^2_{61}$
are $(2,2)$-move equivalent to $9_{49}$.
Here we show additionally that $9_{49}$
and $\bar 9_{49}$ are related by one 5-move, in particular they are
$(2,2)$-move equivalent. The 5-move relation between $9_{49}$
and $\bar 9_{49}$ is illustrated in Figure 2.3.\\
Part (iv) has been proven in \cite{DIP} except the fact that
the closure of $(\sigma_1\sigma_2)^6$ and of $(\sigma_1\sigma_2)^{-6}$
are $(2,2)$-move equivalent. It is the case because, as noted in \cite{DIP},
the closure of $(\sigma_1\sigma_2)^6$ is $(2,2)$-move equivalent to the
knot $9_{49}$.
\end{proof}

\centerline{\psfig{figure=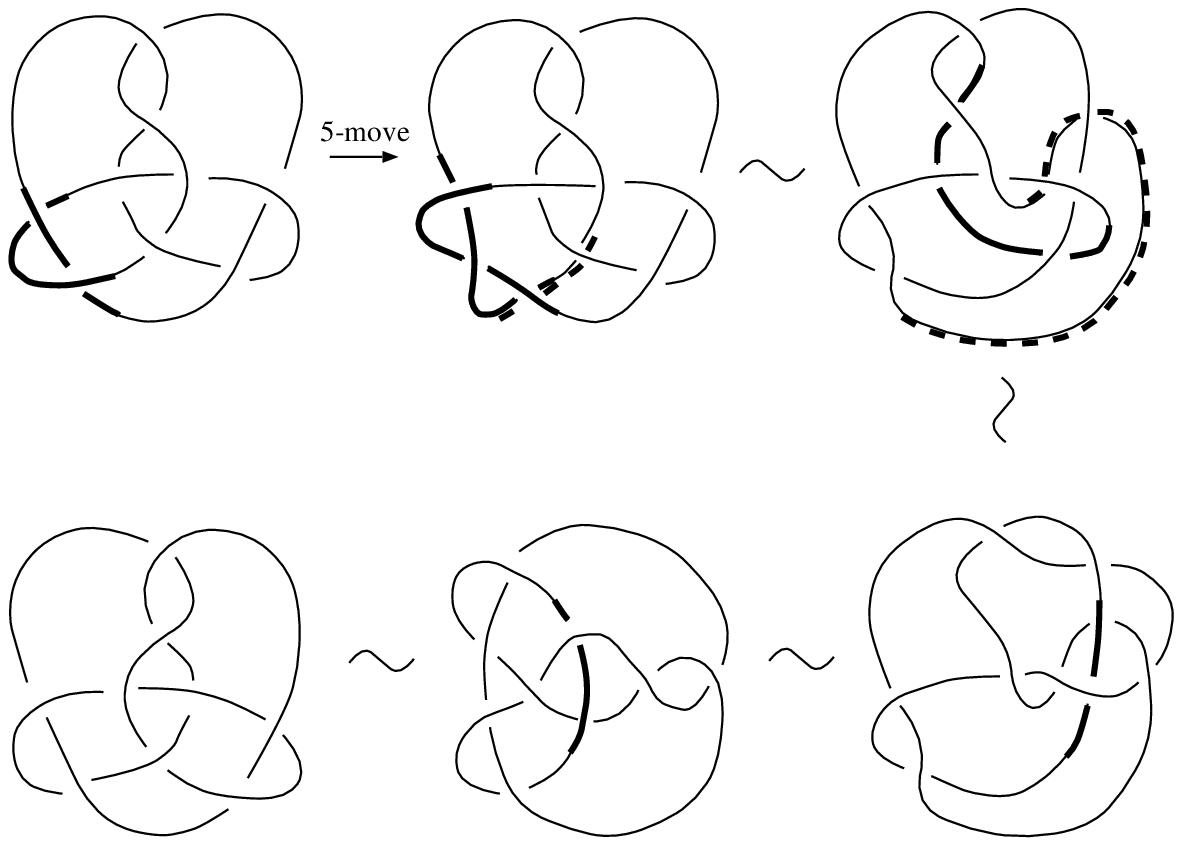,height=6.6cm}}
\centerline{Figure 2.3; $9_{49}$ to $\overline{9_{49}}$}
\ \\

It remains the open problem whether $9_{40}$, $\bar 9_{40}$,
 $9_{49}$ are in different $(2,2)$-move equivalence classes; their
fifth Burnside groups are the same.

\subsection{Kauffman polynomial and $(2,2)$-moves}\label{Subsection 2.4}

It was noted in \cite{Pr-4,Pr-5} that for links $L$ and $L'$ related by one
$(2,2)$-move, and the Kauffman polynomial\footnote{For $a=1$ the Kauffman
polynomial $F_L(a,x)$, was developed before, at the beginning of 1985
by Brandt, Lickorish, Millett and Ho \cite{BLM,Ho}, and denoted by $Q_L(x)$.
It satisfies the skein relation
$Q_{L_+} + Q_{L_-} = x(Q_{L_0}+Q_{L_{\infty}})$, compare Subsection 3.2.}
one has $F_{L'}(1,2cos(2{\pi}/5))=-F_L(1,2cos(2{\pi}/5))$
and that $5(F_L(1,2cos(2{\pi}/5)))^2={col}_5(L)$.
Invariants of $(2,2)$-moves are also invariants of $5$-moves. Furthermore
we can gain some more information from the fact that $(2,2)$-move
is changing the sign of the Kauffman polynomial $F_L(1,2cos(2{\pi}/5))$.
From this it follows that if two links $L$ and $L'$ are $(2,2)$-move
equivalent then the number of moves needed to go from one to another
is even if and only if $F_L=F_{L'}$.
In particular, because $5$-move is a composition of two $\pm(2,2)$-moves,
we have:
\begin{lemma}\label{Lemma 2.6}
(i) If two links differ by an odd number of $\pm(2,2)$-moves, then they
are not $5$-move equivalent. \\
(ii) $F_L(1,2cos(2\pi/5))$ is an invariant of $5$-moves.
\end{lemma}

As a corollary we are able now to prove a variant of Lemma 2.3 for $5$-moves.
\begin{corollary}
The twelve $2$-tangles $[\infty],[0],[-1],[1],[-2],[2],[\frac{2}{5}],
[\frac{5}{2}],[\frac{3}{2}],[-\frac{3}{2}],[\frac{1}{2}]$ and
$[-\frac{1}{2}]$ (Figure 2.4) are in different classes of $5$-move equivalence.
\end{corollary}

\begin{proof}
By Lemma 2.3 the first $6$ tangles in the list are
in different classes of $(2,2)$-move equivalence.
The other six $2$-tangles differ from the first six by a single $(2,2)$-move.
\end{proof}

\centerline{\psfig{figure=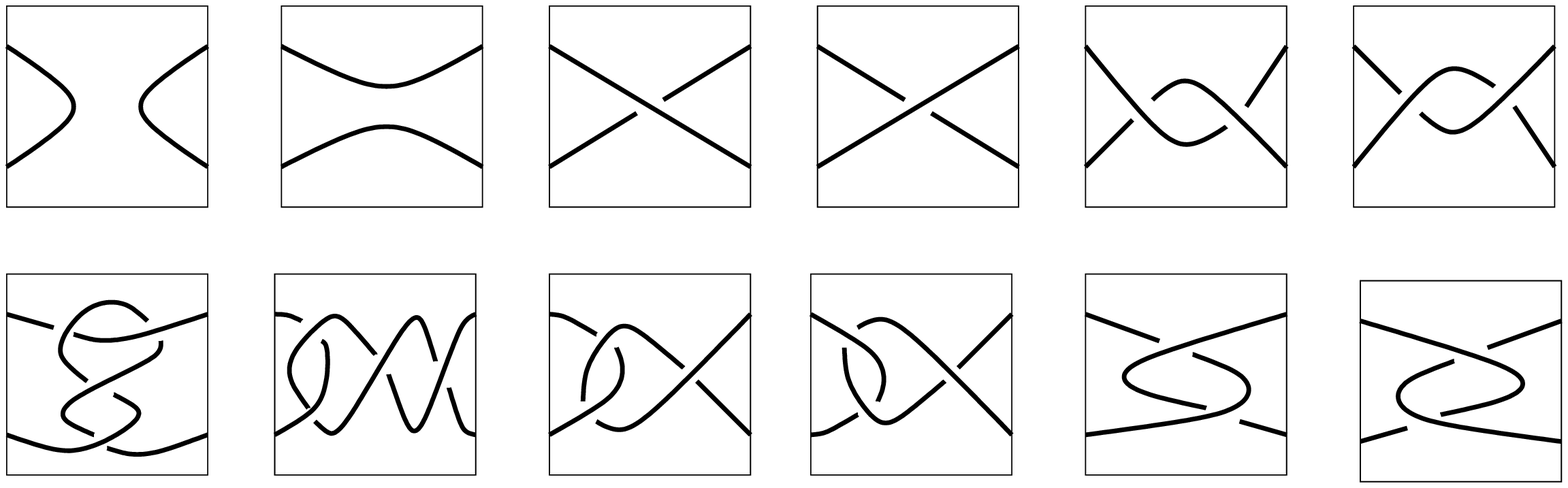,height=4cm}}
\centerline{Figure 2.4;  Basic $12$ tangles}

Lemma 2.6 also allows us to improve slightly the statement of Example 2.4,
 in case of $5$-moves:
\begin{example}\label{Example 2.8}
$L(T_A,k)$ is not $5$-move equivalent to a trivial link if
$T_A$ is (2,2)-move equivalent to crossless $2$-tangle by $n$
$\pm (2,2)$-moves,
and $n+k$ is odd. It is the case because $L(T_A,k)$ can be reduced
to a trivial link by an odd number of $\pm (2,2)$-moves.
\ Notice that the pretzel tangle
$PT_{[m[2],[s]]}= \underbrace{[\frac{1}{2}]*...*[\frac{1}{2}]}_{m\ times}*[s]$
(compare Figure 5.1)
can be changed by $m$ $(2,2)$-moves to $[-2m+s]$-tangle.
\end{example}

We discuss the Kauffman polynomial and $5$-moves in detail, in Section $3.2$.

\subsection{Classification of rational tangles and links
            up to $5$-move equivalence}

In this subsection we classify 5-move equivalence
classes of rational tangles and links (Lemma 2.10 and Theorem 2.9).
From this we get in Section 5 a partial
classification of Montesinos links (including complete classification
of pretzel links) Theorem 5.1.
\begin{theorem}\label{Theorem 2.9}
\begin{enumerate}
\item[(i)]
A rational link is 5-move equivalent to the trivial knot ($T_1$),
the trivial link of $2$ components ($T_2$), the Hopf link ($H$),
or the figure eight knot ($4_1$).
\item[(ii)]
Two rational links are 5-move equivalent if and only if
they have the same value of $F_L=F_L(1,2cos(\frac{2\pi}{5}))$.
The values for $T_1$, $T_2$, $H$, and $4_1$ are: $F_{T_1}=1$,
$F_{T_2}=\sqrt 5$, $F_{H}=-1$, $F_{4_1}=-\sqrt 5$.
\end{enumerate}
\end{theorem}

In Section 5 we will see that rational tangles are also classified
by the absolute value of the Jones polynomial,
$V(L)=|V_L(e^{\frac{\pi i}{5}})|$; compare Table 7.1.

We discuss more of the use of Kauffman polynomial in analysis of 5-move
equivalence in Subsection 3.2.

We deduce Theorem 2.9 from the more general result about
 rational tangles.
\begin{lemma}\label{Lemma 2.10}
Every rational tangle can be reduced by 5-moves to one of twelve 2-tangles
in Figure 2.4 (they are: $[\frac{1}{0}]$, $[0]$,
$[-1]$, $[1]$, $[-2]$, $[2]$, $[\frac{2}{5}]$, $[\frac{5}{2}]$,
$[\frac{3}{2}]$,
$[-\frac{3}{2}]$, $[\frac{1}{2}]$, and $[-\frac{1}{2}]$
 Furthermore these 12 tangles are
representing different 5-move equivalence classes.
\end{lemma}

Before we prove Lemma 2.10 we give an easy to use rule to
recognize quickly to which of 12 tangles the given $\frac{p}{q}$-tangle
is $5$-move reducible.
\begin{proposition}\label{Proposition 2.11} Every rational tangle
$[\frac{p}{q}]$ is in one of twelve $5$-move classes of Lemma 2.10
according to the following rules.
\begin{enumerate}
\item[$(\frac{1}{0})$]\ $q \equiv 0 \mod 5$ and $p \equiv \pm 1 \mod 5$.
\item[$(\frac{2}{5})$]\ $q \equiv 0 \mod 5$ and $p \equiv \pm 2 mod 5$.
\item[$(\frac{0}{1})$]\ $q \equiv \pm 1 \mod 5$ and $p \equiv 0 \mod 5$.
\item[$(\frac{5}{2})$]\ $q \equiv \pm 2 \mod 5$ and $p \equiv 0 \mod 5$.
\item[$(-\frac{1}{1})$]\ $q \equiv \pm 1 \mod 5$ and $p \equiv -q \mod 5$.
\item[$(\frac{3}{2})$]\ $q \equiv \pm 2 \mod 5$ and $p \equiv -q \mod 5$.
\item[$(\frac{1}{1})$]\ $q \equiv \pm 1 \mod 5$ and $p \equiv q \mod 5$.
\item[$(-\frac{3}{2})$]\ $q \equiv \pm 2 \mod 5$ and $p \equiv q \mod 5$.
\item[$(-\frac{2}{1})$]\ $q \equiv \pm 1 \mod 5$ and $p \equiv -2q \mod 5$.
\item[$(\frac{1}{2})$]\ $q \equiv \pm 2 \mod 5$ and $p \equiv -2q \mod 5$.
\item[$(\frac{2}{1})$]\ $q \equiv \pm 1 \mod 5$ and $p \equiv 2q \mod 5$.
\item[$(-\frac{1}{2})$]\ $q \equiv \pm 2 \mod 5$ and $p \equiv 2q \mod 5$.
\end{enumerate}
We can say succinctly that
two rational tangles $[\frac{p}{q}]$ and $[\frac{p'}{q'}]$,
are 5- move equivalent if and only if \\
$q\equiv q' \mod 5$ and $p \equiv p' \mod 5$, or \\
$q\equiv -q' \mod 5$ and $p \equiv -p' \mod 5$.
\end{proposition}
\begin{corollary}\label{Corollary 2.12}
A rational link of type $\frac{p}{q}$ is
\begin{enumerate}
\item[(i)] 5-move equivalent to the trivial link of 2 components iff
$p \equiv 0 \mod 5$ and $q \equiv \pm 1 \mod 5$.
\item[(ii)] 5-move equivalent to the figure eight knot iff
$p \equiv 0 \mod 5$ and $q \equiv \pm 2 \mod 5$.
\item[(iii)] 5-move equivalent to the trivial knot iff
$p \equiv \pm 1 \mod 5$
\item[(iv)] 5-move equivalent to the Hopf link iff
$p \equiv \pm 2 \mod 5$.
\end{enumerate}
\end{corollary}
\begin{proof} It suffices to use Proposition 2.11 when analyzing all
12 tangles of Figure 2.4.  The rational $\frac{p}{q}$ link is  the
numerators  of the tangle $[\frac{p}{q}]$.
\end{proof}

Proposition 2.11 follows from the proof of Lemma 2.10 and in particular
from the fact that in the reduction of any rational tangle to one of 12
tangles we stay in the family of rational tangles and the terms of
related continued fractions are preserved modulo 5.

As a preparation for the proof of Lemma 2.10 we show
\begin{proposition}
The rational tangle $[\frac{3}{2}]$ is 5-move equivalent to the tangle
$[\frac{2}{3}]$, and similarly $[-\frac{3}{2}]$ is 5-move equivalent
to the tangle $[-\frac{2}{3}]$. The rational tangle $[\frac{5}{2}]$
is 5-move equivalent to the tangle $[-\frac{5}{2}]$, and
similarly $[\frac{2}{5}]$ is 5-move equivalent
to the tangle $[-\frac{2}{5}]$. Furthermore,
the rational tangle $[\frac{5}{3}]$ is 5-move equivalent to the tangle
$[\frac{5}{2}]$, and the tangle $[-\frac{5}{3}]$ is 5-move
equivalent to the tangle $[-\frac{5}{2}]$. Similarly,
the rational tangle $[\frac{3}{5}]$ is 5-move equivalent to the tangle
$[\frac{2}{5}]$, and the tangle $[-\frac{3}{5}]$ is 5-move
equivalent to the tangle $[-\frac{2}{5}]$.
\end{proposition}
\begin{proof} The transformation of $[\frac{3}{2}]$ to $[\frac{2}{3}]$
is illustrated in Figure 2.5. Algebraically we have $[\frac{3}{2}] =
[1+\frac{1}{2}] \stackrel{5}\leftrightarrow [1-\frac{1}{3}] =
[\frac{2}{3}]$; we use  $\stackrel{5}\leftrightarrow $ to denote a
transformation by one $\pm 5$-move.
The transformation of $[\frac{5}{2}]$
to $[-\frac{5}{2}]$ is illustrated in Fig. 2.6; algebraically we have
$[\frac{5}{2}] = [2+\frac{1}{2}] \stackrel{5}\leftrightarrow
[-3+\frac{1}{2}] = [-\frac{5}{2}]$.
 Finally, the transformation
of $[\frac{5}{3}]$ to $[\frac{5}{2}]$ is illustrated in Figure 2.7;
algebraically we have $[\frac{5}{3}] = [2-\frac{1}{3}]
\stackrel{5}\leftrightarrow [2+\frac{1}{2}] = [\frac{5}{2}]$.
Other cases of Proposition 2.13 can be obtained from the above by
rotation and mirror images.
\end{proof}
\centerline{\psfig{figure=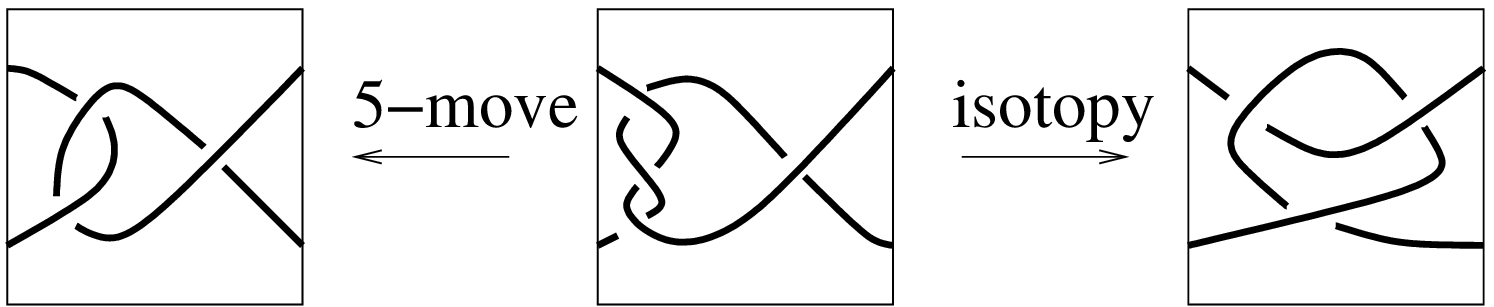,height=1.5cm}}

\begin{center}
Fig. 2.5
\end{center}
\ \\
\centerline{\psfig{figure=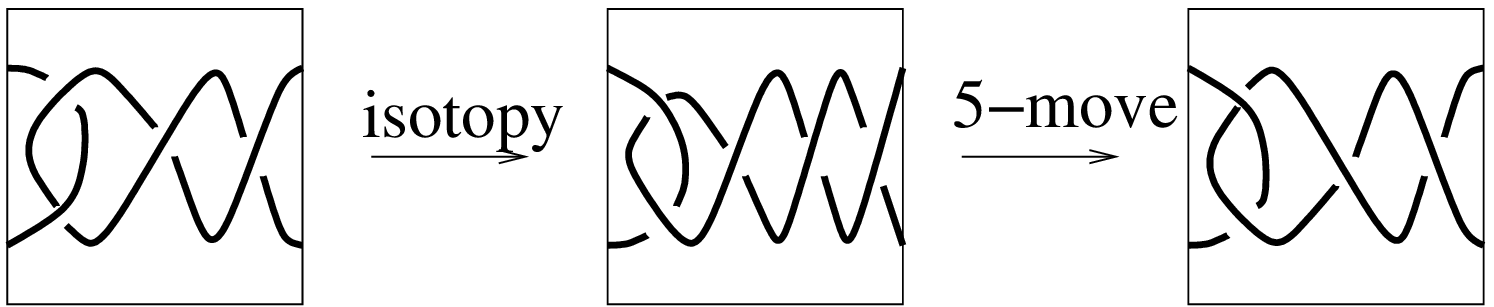,height=1.5cm}}
\begin{center}
Fig. 2.6
\end{center}
\ \\
\centerline{\psfig{figure=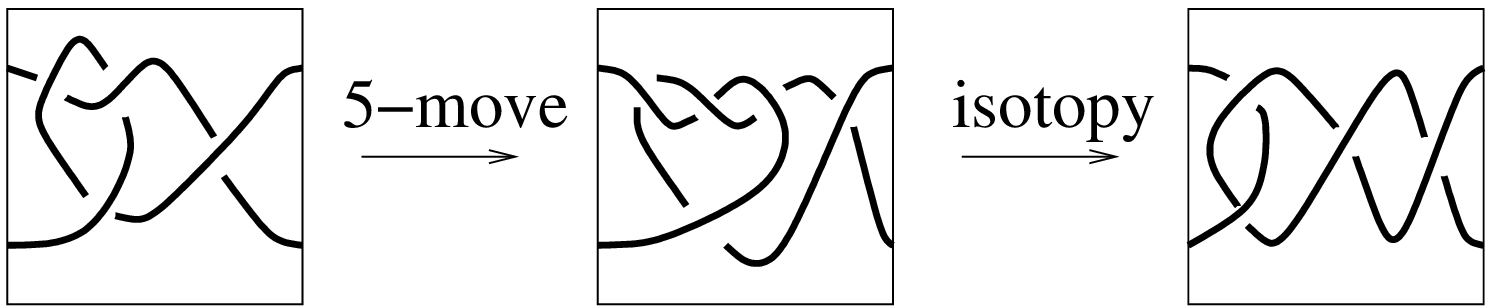,height=1.5cm}}
\begin{center}
Fig. 2.7
\end{center}

\begin{proof} We prove Lemma 2.10 by induction on the minimal number
of crossing of a rational tangle. To make our proof short we use the fact
(version of the Tait conjecture) that the minimal number of crossings
is realized by an alternating diagram (in a continued fractional expansion
it is reflected by a fact that all entries are nonnegative or all
are nonpositive) and that non-alternating
diagram of a rational tangle cannot realize the minimal number of crossings.\\
For diagrams with no more than 4 crossings the result holds by
Proposition 2.13 as any reduced alternating diagram of a rational tangle is
either listed in Lemma 2.10 or in Proposition 2.13.
We assume now that Lemma 2.10 holds for rational tangles of at most
$n$ crossings ($n \geq 4$) and let a rational
tangle $T$ has $n+1$ crossings.
$T$ is obtained from a tangle $T'$ by adding
one crossing. By inductive assumption
 we can reduce $T'$ by 5-moves to one of 12 tangles from lemma 5.2.
Then $T$ is reduced to a tangle $T''$  of at most 5 crossings. If
$T''$ has less than 5 crossings or  is a non-alternating tangle we can
use the fact that lemma is proven already for tangles of up to 4 crossings.
Otherwise, $T'$ was  reduced
to $[\frac{5}{2}]$ or $[\frac{2}{5}]$ tangles and $T''$ is alternating.
We can, however, change by a 5-move the tangle $[\frac{5}{2}]$
(resp. $[\frac{2}{5}]$) to  $[-\frac{5}{2}]$ (resp. $[-\frac{2}{5}]$)
resulting
in non-alternating tangle with 5 crossings which is $5$-move
 reducible to a tangle with no more than 4 crossings for which
Lemma 2.10 already holds.
\end{proof}

\section{Invariants of $5$-moves and their applications}
Invariants of $(2,2)$-moves are also invariants of $5$-moves
and we can employ them as the first step in analyzing links up to $5$-moves.
In this section we use Jones, Kauffman bracket, and Kauffman polynomials
for more detailed analysis of links up to $5$-moves.

\subsection{Jones polynomial and Kauffman bracket of $5$-moves}
In this subsection we use the Jones polynomial and its
Kauffman bracket version to analyze 5-moves.
We work with unoriented diagrams so the Jones
 polynomial $V_L(t)$ is well defined only up to an invertible elements
of $Z[t^{\pm 1/2}]$.
We can do slightly better and define
$\tilde V_L(t) = (t^{\frac{3}{2}})^{-lk(L)}V_L(t)$ which does
not depend on orientation of $L$.
 We use this version of the Jones polynomial
in Section 5.2.

We start from the general formula about the $k$-move and the Kauffman bracket
polynomial and the Jones polynomial. We base our summary on  \cite{Pr-1}.
Recall that the Kauffman bracket polynomial of a link diagram,
$\langle  L \rangle \in Z[A^{\pm 1}]$, satisfies
the Kauffman bracket skein relation \cite{Kau}:
$$\langle  L_+ \rangle = A \langle  L_0 \rangle +
A^{-1}\langle  L_{\infty} \rangle.$$ Let $L_k$
(\parbox{1.6cm}{\psfig{figure=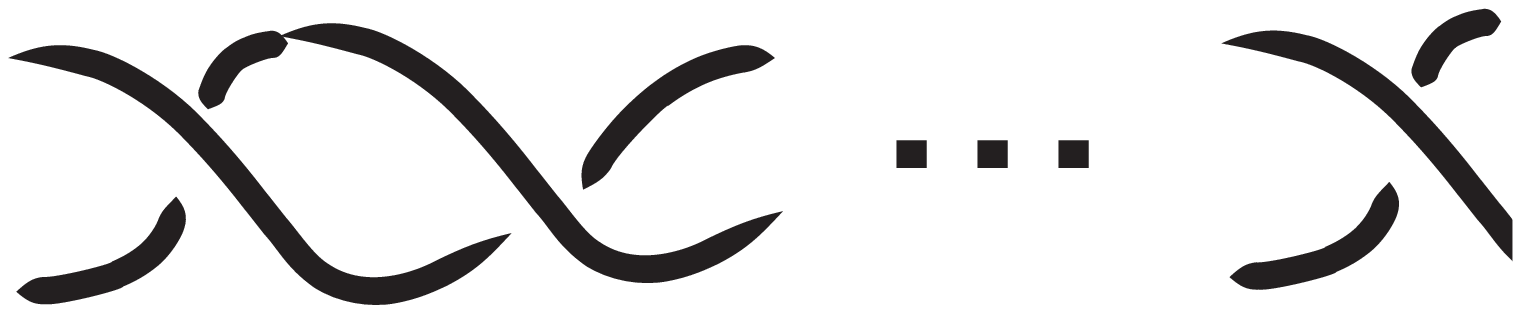,height=0.6cm}}) be obtained from
$L_0$ ($\sg $) by a $k$-move\footnote{We draw the parts of the diagrams which
are involved in the move. The convention for ``$k$-move" used here
 is well rooted in knot theory literature but we should remember that our
$k$-move is a rational $-\frac{k}{1}$-move in Conway's notation.
We also denote
\parbox{0.6cm}{ \psfig{figure=L+nmaly.eps,height=0.5cm}}
by $L_+$ as in Fig. 1.2, but, following Conway, we call
this the $[-1]$ tangle.}
  ($k$ right-handed half-twists added). We have:
\begin{proposition}\label{Proposition 3.1}
\begin{enumerate}
\item[(i)] $\langle  L_k \rangle  =
A^k \langle  L_0\rangle +
A^{-k}\frac{A^{2k} - (-1)^kA^{-2k}}{A^2+A^{-2}}
\langle  L_{\infty}\rangle $.
\item[(ii)] In particular, $\langle  L_5\rangle  =
A^5\langle  L_0\rangle +
A^{-13}(\frac{A^{20}+1}{A^4+1})\langle  L_{\infty}\rangle$, and
$\langle  L_5\rangle  \equiv
A^5\langle  L_0\rangle \mod \frac{A^{20}+1}{A^4+1}$.
\end{enumerate}
\end{proposition}

If we work with framed unoriented links then a $k$-move changes the
Kauffman bracket by $A^5$, modulo $\frac{A^{20}+1}{A^4+1}$. However,
when working with unoriented unframed links then $\langle  L\rangle$ modulo
$\frac{A^{20}+1}{A^4+1}$ is preserved only up to the power of $\pm A^i$.
We write $I_A=\frac{A^{20}+1}{A^4+1}$ and $f(A)\ea g(A)$ if
$f(A) \equiv \pm A^i g(A) \mod I_A$ for some $i$.

The Jones polynomial $V_L(t) \in Z[t^{\pm \frac{1}{2}}]$ can be obtained from
the Kauffman bracket polynomial by putting $t=A^{-4}$ in
$(-A^3)^{-w(L)}\langle  L\rangle$, where $L$ is equipped with any
orientation and $w(L)$ is the writhe or Tait number of an oriented
diagram $L$ ($w(L)={\sum}_p sgn (p)$ where the sum is taken over all
crossings $p$ of oriented diagram $L$. Similarly,
$\tilde V_L(t)= (-A^3)^{-sw(L)}\langle  L\rangle$, for $t=A^{-4}$ and
the self-writhe number $sw(L)$ of an unoriented diagram $L$ is
$sw(L) = w(L) - 2lk(L) = {\sum}_p sgn (p)$ where the sum is taken over all
self-crossings $p$ of unoriented diagram $L$.

\begin{corollary}[\cite{Pr-1}]\label{Corollary 3.2}
\begin{enumerate}
\item[(i)] $V_{L_5}(t) \equiv \pm t^{i/2} V_{L_0} \mod I_t$, for
some $i$, where $I_t = (\frac{t^{5}+1}{t+1})$. We write
succinctly, $V_{L_5}(t) \et V_{L_0} \et \tilde V_{L_5}(t)$.
\item[(ii)] For $t=e^{\pi i/5}$, $|V_L(t)|=|\tilde V_L(t)|$ is
an invariant of $5$-move equivalence classes of links.
We denote this invariant by $V(L)$; compare Tables 4.1 and 7.1.
\item[(iii)] The space $Z[t]/(t^4=t^3-t^2+t-1)$ is isomorphic to the space
of polynomials of degree at most 3. The Jones polynomial $V_L(t)$ is
either in $Z[t^{\pm1}]$ (if $L$ has odd number of components) or
$t^{1/2}V_L(t)$ is in $Z[t^{\pm 1}]$ (if $L$ has even number of components).
If one reduces this polynomial modulo $I_t$ and then takes the
result  up to $\pm t^i$ one gets the set of $5$ polynomials
(up to the sign). We denote this set by $V(L,5)$. Then
if two links are 5-move equivalent then they have
the same (up to the sign) set of polynomials $V_L(t,5)$.
\end{enumerate}
\end{corollary}
\begin{example}\label{Example 3.3}\ \\
(i) $V(L)$ classifies rational links. We have $V(T_1)=1$,
$V(T_2)= 2cos(\pi/10) \approx 1.90211$,
$V(H) = 2cos(\pi/5) \approx 1.61803$, $V(4_1)=0$. In particular,
the $\frac{p}{q}$-rational link is $5$-move equivalent to $H$ iff
$p\equiv \pm 2 \mod 5$; compare Corollary 2.12.\\
(ii) For the pretzel link $6^3_1$ and its mirror image $\bar 6^3_1$
we have $V(6^3_1)=V(\bar 6^3_1) \approx 2.497$ but
$V_{6^3_1}(t,5)=\{2+t^2,2t+t^3,-1+t+t^2+t^3,-1+2t^3, -2+t-2t^2+2t^3\}
\neq V_{\bar 6^3_1}(t,5)=\{1+2t^2,t+2t^3,-2+2t-t^2+2t^3,-2+t^3,-1-t-t^2+t^3\}$.
In particular, $V(6^3_1)$ is neither $5$-move equivalent to its mirror
image nor to any rational knot; compare Example 5.14 and Proposition 5.16.
\end{example}

\subsection{Using Kauffman polynomial to analyze 5-moves}
Recall that the 2-variable Kauffman polynomial of regular
isotopy of link diagrams (or equivalently of framed links with
blackboard framing) $\Lambda_L(a,x) \in Z[a^{\pm},x^{\pm}]$
is defined recursively as follows \cite{Kau}:\\
(i) (Initial condition) $\Lambda_{\bigcirc}(a,x)=1$.
\\
(ii) (First Reidemeister move, or framing condition)
$\Lambda_{\ \parbox{0.4cm}{ \psfig{figure=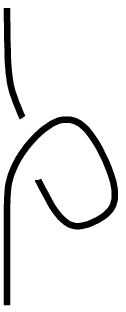,height=0.6cm}}}
(a,x)=
a \Lambda_{\ \parbox{0.2cm}{\psfig{figure=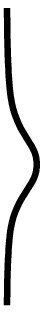,height=0.6cm}}}
(a,x)$.
\\
(iii) (Kauffman skein relation) $ {\Lambda}_{L_{+}}(a,x) +
\Lambda_{L_{-}}(a,x)=x(\Lambda _{L_{0}}(a,x)+\Lambda _{L_{\infty }}(a,x))$.
\\
The Kauffman polynomial $F_L(a,x)$ of oriented links,
 is obtained by normalizing ${\Lambda}_L(a,x)$, that is
$F_L(a,x)=a^{-w(L)}{\Lambda}_{L}(a,x)$,
 where $w(L)={\sum}_p sgn p$ where the sum is taken over all crossings $p$ of
oriented diagram $L$.

Let $L_k$ be a diagram obtained from $L=L_0$ by a $k$-move.
In \cite{Pr-1} we derived the following formula.
\begin{theorem}\cite{Pr-1}\label{Theorem 3.3}
$$\Lambda _{L_k}(a,x) = v_1^{(k)}(x)\Lambda_{L_1}(a,x)
- v_1^{(k-1)}(x)\Lambda_{L_0}(a,x)+xv_2^{(k)}(a,x)\Lambda_{L_{\infty}}(a,x),$$
where the polynomials
$v_1^{(k)}(x)$ and $v_1^{(k-1)}(x)$   are (shifted) Chebyshev polynomials of
the first type{\footnote{Chebyshev polynomial of the first type
$T_k(x)$ satisfies:
$T_k(x)=T_{k-1}(x)-T_{k-2}(x), T_0(x)=1,T_1(x)=x$ and for $x=p+p^{-1}$ we
have $T_k(x)= \frac{p^{k+1}-p^{-k-1}}{p-p^{-1}}$.
Therefore $v_1^{(k)}(x)=T_{k-1}(x)$ and $v_1^{(k-1)}(x)=T_{k-2}(x)$.
Furthermore, $v_2^{(k)}(a,x)$ is a generating functions of Chebyshev
polynomials, that is $v_2^{(k)}(a,x)={\sum}_{i=1}^{k-1}T_{i-1}(x)a^{i-k}$
(see Corollary 3.5).}}.
By putting $x=p+p^{-1}$ we obtain
$v_1^{(k)}(x)=\frac{p^k-p^{-k}}{p-p^{-1}}$,
$v_1^{(k-1)}(x)=\frac{p^{k-1}-p^{1-k}}{p-p^{-1}} = pv_1^{(k)}(x) - p^k$,
and $v_2^{(k)}(a,x)=\frac{(-a^{-1}(p^k-p^{-k})+
p(a^{-k}-p^{-k})-p^{-1}(a^{-k}-p^k))}{(p-p^{-1})(a+a^{-1}-(p+p^{-1}))}$.
\end{theorem}

Theorem 3.4 can be reformulated in the following, useful for our analysis,
form.
\begin{corollary}\label{Corollary 3.5} \
For $x=p+p^{-1}$ we have
$$\Lambda _{L_k}(a,x) = \frac{p^k-p^{-k}}{p-p^{-1}}\Lambda_{L_1}(a,x) -
\frac{p^{k-1}-p^{1-k}}{p-p^{-1}}\Lambda_{L_0}(a,x) +$$
$$xa^{-1}\bigl({\sum}_{i=0}^{k-2} (\frac{p^{i+1}-p^{-i-1}}{p-p^{-1}})a^{i-k+2}
\bigr)
\Lambda_{L_{\infty}}(a,x).$$
Furthermore, the coefficient of
the last summand reduces, for $a=p$ (so also $x=a+a^{-1})$, to
$$xa^{-1}({\sum}_{i=0}^{k-2} \frac{a^{i+1}-a^{-i-1}}{a-a^{-1}}a^{i-k+2}) =
xa^{k-3}({\sum}_{i=1}^{k-1} i(a^{-2})^{i-1}) =
xa^{k-3}\frac{d(\frac{(a^{-2})^k-1}{a^{-2}-1})}{d(a^{-2})}.$$
\end{corollary}
\begin{proof}
Corollary 3.5 can be derived directly from Theorem 3.4 but it can be also
quickly proven by induction on $k$. The inductive step has the form:
$$\Lambda_{L_{k+1}}(a,p+p^{-1}) = (p+p^{-1})\Lambda_{L_k}(a,p+p^{-1}) -
\Lambda_{L_{k-1}}(a,p+p^{-1})+
a^{-k}(p+p^{-1})\Lambda_{L_{\infty}}(a,p+p^{-1})$$
$$=((p+p^{-1})(\frac{p^{k}-p^{-k}}{p-p^{-1}})-
(\frac{p^{k-1}-p^{1-k}}{p-p^{-1}})) \Lambda_{L_1} -
((p+p^{-1})(\frac{p^{k-1}-p^{1-k}}{p-p^{-1}})-
(\frac{p^{k-2}-p^{2-k}}{p-p^{-1}})) \Lambda_{L_0}+ $$
$$\bigl(xa^{-1} \bigl((p+p^{-1}){\sum}_{i=0}^{k-2}
(\frac{p^{i+1}-p^{-i-1}}{p-p^{-1}})a^{i-k+2})-
{\sum}_{i=0}^{k-3} (\frac{p^{i+1}-p^{-i-1}}{p-p^{-1}})a^{i-k+3}\bigr)
+ a^{-k}(p+p^{-1})\bigr) \Lambda_{L_{\infty}} $$
$$=\frac{p^{k+1}-p^{-k-1}}{p-p^{-1}}\Lambda_{L_1} -
\frac{p^{k}-p^{-k}}{p-p^{-1}}\Lambda_{L_0}+
xa^{-1}({\sum}_{i=0}^{k-1} (\frac{p^{i+1}-p^{-i-1}}{p-p^{-1}})a^{i-k+1})
\Lambda_{L_{\infty}}
$$
\end{proof}

Let the ideal $I_{v_1^{(k)},v_2^{(k)}}= (v_1^{(k)}(x),v_2^{(k)}(a,x))$
be the ideal in $Z[a^{\pm},x^{\pm}]$ generated by $v_1^{(k)}(x)$ and
$v_2^{(k)}(a,x)$.
\begin{corollary}\label{Corollary 3.6} \ \ \ \
(i) If two framed links $L$ and $L'$
differ by a $k$-move then
 ${\Lambda}_L(a,x) \equiv a^k {\Lambda}_{L'}(a,x) \mod
I_{v_1^{(k)},v_2^{(k)}}$; additionally $a^{2k} \equiv 1 \mod
I_{v_1^{(k)},v_2^{(k)}}$.\\
(ii) $p^k \equiv a^k \mod I_{v_1^{(k)},v_2^{(k)}}$ and
if $a+a^{-1} -x$ is invertible in our ring (e.g. we consider
$Z[a^{\pm 1}, x^{\pm 1}, (a+a^{-1} -x)^{-1}]$) then
$I_{v_1^{(k)},v_2^{(k)}}= (p^{2k}-1, p^k-a^k)$. In particular,
for any numbers
$a_0, p_0 \in \C$, $a_0^{2k}=1, p_0^k=a_0^k$, $a_0\neq p_0,p_0^{-1}$,
$p_0 \neq 1, -1, i, -i$, we have ${\Lambda}_L(a_0,x_0) =
(a_0)^k {\Lambda}_{L'}(a_0,x_0)$.\\
(iii) If $a_0=1, p_0^k=1$, $p_0 \neq 1, -1, i, -i$ then
${\Lambda}_L(a_0,x_0)$ is a $k$-move equivalence invariant of
unoriented, unframed links. \\
(iv) In the case of $x=a+a^{-1}$, the ideal $I_{v_1^{(k)},v_2^{(k)}}\subset
Z[a^{\pm},x^{\pm}] $ reduces to $I_{x=a+a^{-1}}\subset Z[a^{\pm}]$,
where $I_{x=a+a^{-1}} = (v_1^{(k)}(a),v_2^{(k)}(a,a+a^{-1})=
(\frac{(a^{-2})^k -1}{a^{-2}-1},
\frac{d(\frac{(a^{-2})^k-1}{a^{-2}-1})}{d(a^{-2})})=
(k,1+2a^{-2}+3a^{-4} +...+ (k-1)a^{4-2k}) = (k,1+2a^2+3a^4+...+
(k-1)a^{2k-4}).$ Furthermore, for $k$ a prime number,
$I_{x=a+a^{-1}} = (k, (a^2-1)^{k-2})$.
\end{corollary}
\begin{proof} Corollary 3.6(i)-(iii) is proven in \cite{Pr-1}; here let us
only notice that
$(a+a^{-1}-x)v_2^{(k)}(a,x)= (p-a^{-1})v_1^{(k)}(x) + a^{-k}-p^k$ which
allows short proof of (ii) and (iii) from (i). Furthermore, notice that
$F_L(a,x) = F_L(-a,-x)$ and, equivalently,
 ${\Lambda}_L(a,x)=(-1)^{w(L)}\Lambda_L(-a,-x)$. This, for odd $k$ allows
us to consider only substitution $a_0^k=1=p_0^k$ in Corollary 3.6(ii).
To prove (iv), we use the second part of Corollary 3.5, where
it was noted that $v^{(2)}_2(a,a+a^{-1}) =
a^{k-3}({\sum}_{i=1}^{k-1} i(a^{-2})^{i-1}) =
a^{k-3}\frac{d(\frac{(a^{-2})^k-1}{a^{-2}-1})}{d(a^{-2})}=
a^{1-k}(k-1 + (k-2)a^2 + ... +2a^{2k-6}+ a^{2k-4})$. Furthermore, we have
$$k(a^{-2})^{k-1} = \frac{d(a^{-2})^k-1)}{d(a^{-2})}=
\frac{d(\frac{(a^{-2})^k-1}{a^{-2}-1})(a^{-2}-1)}{d(a^{-2})}=
\frac{(a^{-2})^k-1}{a^{-2}-1} +
(a^{-2}-1)\frac{d(\frac{(a^{-2})^k-1}{a^{-2}-1})}{d(a^{-2})}.$$
Furthermore, for $k$ being a prime number $\frac{(a^{-2})^k-1}{a^{-2}-1}
\equiv (a^{-2} - 1)^{k-1}\mod k$ and, therefore,
$\frac{d(\frac{(a^{-2})^k-1}{a^{-2}-1})}{da^{-2}}= (k-1)(a^{-2} - 1)^{k-2}
\equiv -(a^{-2} - 1)^{k-2} \mod k$.

\end{proof}

To analyze $5$-move equivalence of links we are interested
in the case of $k=5$. Then we have:
\begin{corollary}\label{Corollary 3.7}\ \\
 (i)
$\Lambda_{L_5}(a,x) \equiv p^5 \Lambda_{L_0} +
xv_2^{(5)}(a,x)\Lambda_{L_{\infty}} \mod (v_1^{(5)}(x))$, \\
\ \ \  $\Lambda_{L_5}(a,x) \equiv p^5 \Lambda_{L_0} \mod (v_1^{(5)}(x),
v_2^{(5)}(a,x))$.\\
(ii) If $a_0,p_0\in \C$ and $a_0^{10}=1$, $p_0^5=a_0^5$,
$a_0\neq p_0,p_0^{-1}$, $p_0\neq 1,-1,i,-i$,
then $\Lambda_L(a_0,x_0)$ is a $5$-move equivalence invariant of
unoriented, unframed links.\\
(iii) $F_L(1,x)$ modulo the ideal $(x^2-x+1)$ (or equivalently
$\frac{p^5-1}{p-1}$) is a $5$-move equivalence invariant of
unoriented, unframed links. In particular,  if $p_0 = e^{2\pi i/5}$
(that is, $x_0 = 2cos(2\pi/5)$) it is the invariant used in Lemma 2.6(ii).\\
(iv) $\Lambda_{L_5}(a,a+a^{-1}) \equiv a^5\Lambda_{L_0}(a,a+a^{-1})
\mod (5, (a^2-1)^3)$.
\end{corollary}

We are mostly interested, in this paper, in unoriented,
unframed links, so we modify Corollary 3.7 accordingly,
taking onto account the fact that $\Lambda_{L^{(1)}}(a,x)=a\Lambda_L(a,x)$,
where $L^{(1)}$ is a framed link obtained from a framed link $L$ by a
positive twist on the framing of $L^{(1)}$. Let $I_{(a,x)}$ be
an ideal in $Z[a^{\pm 1}, x^{\pm 1}]$ generated by $v_1^{(5)}(x)$ and
$v_2^{(5)}(a,x)$. We write $p(a,x)
\eax
q(a,x)$ if $p(a,x) \equiv a^i q(x) \mod I_{(a,x)}$ for some $i$. Then we have.
\begin{corollary}\label{Corollary 3.8}\ \\
(i) If two links, $L_1$ and $L_2$ are $5$-move equivalent then
$\Lambda_{L_1}(a,x) \eax \Lambda_{L_2}(a,x)$ (equivalently,
$F_{L_1}(a,x) \eax F_{L_2}(a,x)$ as the equality does not depend on
orientation of $L$). In particular:\\
(ii) If $a_0,p_0\in \C$ and $a_0^{10}=1$, $p_0^5=a_0^5$,$a_0\neq p_0,p_0^{-1}$,
 $p_0\neq 1,-1,i,-i$,
then the set $\{a_0^iF_L(a_0,x_0)\}$ is an invariant of $5$-move equivalence
of unframed links. We denote this invariant by $Set(F_L(a_0,x_0))$.
The absolute value of $F_L(a_0,x_0)$ is also a $5$-move invariant. \\
(iii) $F_{L_1}(a,a+a^{-1}) \eaa F_{L_2}(a,a+a^{-1}).  $
\end{corollary}

Let us remark here that if $Set(F_L(a_0,x_0))$ contains the real value, as
in the case of an amphicheiral link, then this value is a 5-move invariant
for $a_0^5=1$ and it is an invariant up to a sign if $a_0^5=-1$.

\begin{example}\label{Example 3.9}
The pretzel link $P_{[2,2,2]}$ is not 5-move equivalent to its
mirror image. We prove it by computing $a^iF_{P_{[2,2,2]}}(a,x)$ for
$a=e^{4\pi i/5}$ and $x=2cos(2\pi/5)$, and checking that it
is never a real number.
\end{example}

Consider links $L_1=4_1\#L'_1$ and $L_2=4_1\#L'_2$, then we have $V_{L_i}(t)=
V_{4_1}(t)V_{L'_i}(t) \equiv 0 \mod I_t$. We can, however, use the
Kauffman polynomial criteria to differentiate, in some cases, $L_1$ from $L_2$.
\begin{example}\label{Example 3.10}\ \\
 The links $4_1\#4_1$ and $4_1\#T_2$ are not 5-move equivalent.
We have $F_{4_1\#4_1}(1,2cos(2\pi/5))= 5$
but $F_{4_1\#T_2}(1,2cos(2\pi/5))=-5$.
\end{example}

\begin{remark}\label{Remark 3.11}
It is an open problem whether the links $L_1=4_1\#4_1\#4_1$ and
$L_2=4_1\#T_2\# T_2$ are  5-move equivalent.  We have: \\
(i) $L_1$ and $L_2$ are $(2,2)$-move equivalent
by two $(2,2)$-moves and
$F_{L_1}(1,2cos(2\pi/5))=-5\sqrt 5 = F_{L_2}(1,2cos(2\pi/5))$. \\
(ii) The criterion of Corollary 3.8(ii) would not separate
$L_1$ and $L_2$ because
if we assume $a_0\neq \pm 1$ then $F_{4_1}(a_0,x_0)=0$.
The last equality follows from the following computation:\\
$F_{4_1}(a,x)= -a^{-2} -1 -a^2 + x(-a^{-1}-a) + x^2(a^{-2}+2+a^2) +
x^3(a^{-1}+a) \stackrel{y=a+a^{-1}}{=} 1-y^2 + xy(x^2+xy-1)= 
1+y(x+y)(x^2-1)$. Then
$(x-y)F_{4_1}(a,x) =
(x^2+x-1)\bigl(y(x^2-x+1)-y^3\bigr) + (y^2+y-1)(xy-x) \equiv 0
\mod(x^2+x-1,y^2+y-1)$. For $x= p+p^{-1}$, $y=a+a^{-1}$ we have
$x^2+x-1= p^{-2}(\frac{p^5-1}{p-1})$ and $y^2+y-1= a^{-2}(\frac{a^5-1}{a-1})$.
Furthermore, $F_{4_1}(a,x)=F_{4_1}(-a,-x)$ and if $a_0^5=p_0^5=1$ then
$(-a_0^5)=(-p_0^5)=-1$, thus for any substitutions from Corollary 3.8(ii),
$F_{4_1}(a_0,x_0)=0$ as long as $a_0\neq \pm 1$.\\
(iii) The criterion of Corollary 3.8(iii) would not separate
$L_1$ and $L_2$ because $F_{T_2}(a,a+a^{-1})= 0= F_{L_2}(a,a+a^{-1})$ and,
less obviously, $F_{L_1}(a,a+a^{-1}) \equiv 0 \mod (5,(a^2-1)^3)$. To see
the last congruence, we notice that $F_{4_1}(a,a+a^{-1}) =
1- 2(a+a^{-1})^2 + 2(a+a^{-1})^4 = a^{-4}(2+6a^2+9a^4+6a^6+2a^8) \equiv
2a^{-4}(a^4+1)(a^2-1)^2 \mod 5.$\\
We checked generally using the Gr{\"o}bner basis method that
$F_{L_1}(a,x) - F_{L_2}(a,x)$ is in the ideal $I_{a,x}$ thus
the method of Corollary 3.8(i) would not distinguish 5-move equivalence
classes of these links.

\end{remark}

\section{Classification of $3$-braid links up to $5$-move
equivalence}\label{Section 4}

The invariants of $5$ moves introduced in previous sections allows
us to classify 3-braid links up to 5-moves almost completely.
There are at least 23 classes of 5-move equivalence and no more than 25.
We use the names of links from Rolfsen book \cite{Rol} for knots up to
10 crossings and links up to 9 crossings. For links of 10 or 11 crossings
we use Knot-Plot tables \cite{Bar} and for links of 12 or 13 crossings we
use names from Thisthtlethwaite tables (for example $12n_{1958}$ denotes
a non-alternating link of 12 crossings which is $1958$th in \cite{Thi}).
\begin{theorem}\label{Theorem 4.1}\
\begin{enumerate}
\item[(i)]
Every link represented by a closed $3$-braid is $5$-move equivalent
to one of the following 25 links:\ \
$T_1,T_2,T_3,H,H \cup T_1,H\#H$, $4_1,6_1^3,
{\bar{6}}_1^3,6_2^3,6_3^3,{\bar{6}}_3^3, 7_5^2,\\
 8_{18}, 8_{10}^3,
{\bar{8}}_{10}^3,$
$ 8_9^3, 9_{40}^2,{\bar{9}}_{40}^2, 9_{41}^2,{\bar{9}}_{41}^2$,
$L10a163$, $\bar L10a163$, and
$11a_{177}$ and $\bar{11}a_{177}$ (represented by 3-braids
$\sigma_1^2\sigma_2^{-2}
  \sigma_1\sigma_2^{-1}\sigma_1\sigma_2^{-1}\sigma_1\sigma_2^{-2}$ and
$\sigma_1^{-2}\sigma_2^{2}
  \sigma_1^{-1}\sigma_2\sigma_1^{-1}\sigma_2\sigma_1^{-1}\sigma_2^{2}$).
\item[(ii)]
These links represent different $5$-move equivalence classes
with the possible exception of two pairs $8_{10}^3$, ${\bar{8}}_{10}^3$,
and $9^2_{40}$, ${\bar{9}}^2_{40}$.
\end{enumerate}
\end{theorem}

\begin{proof}
We use the Coxeter result that $B_3/({\sigma}_i)^5$ is a finite group
(replacing ${\sigma}_i^5$ by ${\sigma}_i^0$ can be achieved
by a $5$-move) \cite{Cox}.
The quotient group has $45$ conjugacy classes.
We list them all in Table 4.1. For each class (generated by GAP) we
choose a representative which is as short as we are able to find
(we did not prove that they are the shortest).
We provide also the value of invariants of $5$-move equivalence
$F=F_L(1,2cos 2 \pi /5)$, $V(L)$, $V_L(t,5)$ for closures of these braids.

A look at the table shows that each conjugacy class of
$B_3/({\sigma}_i^5) $ is $5$-move equivalent to one of $25$ links
of the theorem. Furthermore, each pair of links with possible
exception of (33) (representing ${\bar{8}_{10}^3}$) and
(35) (representing ${8_{10}^3}$), and (36) (representing $9_{40}^2$) and
(40) (representing ${\bar{9}}_{40}^2$) are separated by listed invariants.
\end{proof}

It is an open problem whether ${8_{10}^3}$ and ${\bar{8}_{10}^3}$ are
5-move equivalent. Similarly 5-move equivalence of $9_{40}^2$ and
${\bar{9}}_{40}^2$ is not yet decided.

In the sixth column of the table, we identify the link being
the closure of a representative and its $5$-move reduction, if any.
For example, (29) has a representative
${{\sigma}_2^{-2}}{\sigma}_1{{\sigma}_2}^{-2}{\sigma}_1^2$ which
represents the link ${\bar{7_5^2}}$. This link is 5-move equivalent
to its mirror image ${7_5^2}$, representing (28). This link, in turn,
can be changed by one $(-5)$-move to the pretzel link
$P_{[2,2,2,1]}$ (i.e. $7^3_1$), see example 5.14(ii)). One more $(-5)$-move
changes this link to its mirror image $P_{[2,2,2,-4]}$ ambient isotopic to
$P_{[-2,-2,-2,-1]}$.\
Similarly (21) has a representative
${\sigma}_2^{-2}\sigma_1^2{\sigma}_2^{-1}{\sigma}_1$
describing the rational knot $6_3$ which is $5$-move equivalent
to the Hopf link $H$.

In the last column we list some interesting representatives of
conjugacy classes in $B_3/(\sigma_i^5)$ different from that
listed in the second column. We pay special attention to powers
of $(\sigma_1\sigma_2)$. In our notation
$L_1 \stackrel{5}\sim L_2$ means
5-move equivalence of links and $L_1 \stackrel{5}\approx L_2$ means the
same conjugacy class in $B_3/(\sigma_i^5)$ and is used only for closed
3-braids.

We end this section with one more question: all closed braids in Table 4.1
have a representative with 10 or less crossings except the pair (43) and
(44) with 11 crossings. Is it possible to reduce these closed braids
to links with $10$ crossings? We know that they are not $5$-move
equivalent to any link of $9$ or less crossings as the only links
which share with them $V(L)$ are 3-component links
 $9^3_{21}$ and its mirror image $\bar 9^3_{21}$ which are algebraic links.
We know that (43) and (44) are separated from algebraic links
(even up to $(2,2)$-move equivalence) by $5$th Burnside group (see
Subsections 2.2 and 2.3).
\ \\
\newpage
\textbf{TABLE 4.1:} \textbf{LIST OF 45 CONJUGACY CLASSES OF }$\mathbf{B}_{3}
\mathbf{/(\sigma }_{1}^{5}\mathbf{)}$\\
{\tiny{
\begin{tabular}{|l|l|l|l|l|l|l|}
\hline
${GAP~CC\#}$ & {braids(shortest)} & ${\tiny{{\bf{F}}}}$& ${V(L)}$ &
${V}_{L}{(t,5)}$ {rep}${.}$ & {Link} in $B_3/(a^5)$ & {interesting rep.}
\\ \hline &&&&&&\\
${(1)}$ & ${Id}$ & 5 &${3.61803 }$ & ${  \{1+2t+t}%
^{2}{,... \}}$ & ${  T}_{3}$ &  \\ \hline &&&&&&\\
${  (2)}$ & ${  \sigma }_{2}$ & $\sqrt{5}$ & ${  1.90211 }$ & ${  %
\{1+t}^{2}{,... \}}$ & ${  T}_{2}$ &  \\ \hline &&&&&&\\
${  (3)}$ & ${  \sigma }_{2}^{-1}$ & $\sqrt{5}$ & ${  1.90211 }$ & $%
{  \{1+t,... \}}$ & ${  T}_{2}$ &  \\ \hline &&&&&&\\
${  (4)}$ & ${  \sigma }_{2}^{2}$ & $-\sqrt{5}$  & ${  3.07768 }$ & $%
{  \{1+t+t}^{2}{  +t}^{3}{,... \}}$ & ${  H\sqcup T}_{1}$ &
 \\ \hline &&&&&&\\
${  (5)}$ & ${  \sigma }_{2}^{-2}$ & $-\sqrt{5}$  & ${  3.07768 }$ & $%
{  \{1+t+t}^{2}{  +t}^{3}{,... \}}$ & ${  H\sqcup T}_{1}$ &
 \\ \hline &&&&&&\\
${  (6)}$ & ${  \sigma }_{2}{  \sigma }_{1}$ & ${  1}$ & $%
{  1}$ & ${  \{1,...\}}$ & ${  T}_{1}$ & ${  (\sigma }_{1}{  %
\sigma }_{2}{  )}^{11}$ \\ \hline &&&&&&\\
${  (7)}$ & ${  \sigma }_{2}^{-1}{  \sigma }_{1}$ & ${  1}$ &
${  1}$ & ${  \{1,...\}}$ & ${  T}_{1}$ & \\ \hline &&&&&&\\
${  (8)}$ & ${  \sigma }_{2}^{2}{  \sigma }_{1}$ & $-1$ & ${  %
1.61803 }$ & ${  \{1+t}^{2}{,... \}}$ & ${  H}$ &  \\ \hline &&&&&&\\
${  (9)}$ & ${  \sigma }_{2}^{-2}{  \sigma }_{1}$ & $-1$ & ${  %
1.61803 }$ & ${  \{1+t}^{2}{,... \}}$ & ${  H}$ & ${  %
(\sigma }_{1}{  \sigma }_{2}{  )}^{-8}$ \\ \hline &&&&&&\\
${  (10)}$ & ${  \sigma }_{2}^{-1}{  \sigma }_{1}^{-1}$ & ${  1%
}$ & ${  1}$ & ${  \{1,...\}}$ & ${  T}_{1}$ & ${  (\sigma }_{1}%
{  \sigma }_{2}{  )}^{-11}$ \\ \hline &&&&&&\\
${  (11)}$ & ${  \sigma }_{2}^{2}{  \sigma }_{1}^{-1}$ & $-1$ & $%
{  1.61803 }$ & ${  \{1+t^{2},... \}}$ & $%
{  H}$ & ${  (\sigma }_{1}{  \sigma }_{2}{  )}^{8}$ \\ \hline &&&&&&\\
${  (12)}$ & ${  \sigma }_{2}^{-2}{  \sigma }_{1}^{-1}$ & $-1$ & $%
{  1.61803 }$ & ${  \{1+t}^{2}{,... \}}$ & ${  H}$ &  \\ \hline &&&&&&\\
${  (13)}$ & ${  \sigma }_{2}^{2}{  \sigma }_{1}^{2}$ & $-1$ & $%
{  2.61803 }$ & ${\{1+t+t}^{2},...\}$ & ${  H\#H}$ &  \\ \hline &&&&&&\\
${  (14)}$ & ${  \sigma }_{2}^{-2}{  \sigma }_{1}^{2}$ & 1 & $%
{2.61803 }$ & ${\{1+t+t}^{2},... \}$ & ${  H\#H}$ &
${  (\sigma }_{1}{  \sigma }_{2}%
{  )}^{5}$ \\ \hline &&&&&&\\
${  (15)}$ & ${  \sigma }_{2}^{-2}{  \sigma }_{1}^{-2}$ &  1& $%
{  2.61803 }$ & ${  \{1+t+t}^{2},... \}$ & ${  H\#H}$ &
\\ \hline &&&&&&\\
${  (16)}=({\overline{20}}$) & ${  \sigma }_{2}^2{  \sigma }_{1}{  \sigma }_{2}^{-1}{  \sigma }_{1}$ & ${1}$ &
${  1%
}$ & ${  \{1,...\}}$ & ${\overline{4}}_{1}^{2}{  \stackrel{5}{\sim} T}_{1}$ &  \\ \hline &&&&&&\\
${  (17)}$ & $({  \sigma }_{1}{  \sigma }_{2}^{-1})^2$ & $-\sqrt{5}$ & ${  0}$ & $%
{  \{0,...\}}$ & ${  4}_{1}$ & \\ \hline &&&&&&\\
${  (18)=({\overline{19}})}$ & ${\sigma}_{1}{  \sigma }_{2}^{-2}{  \sigma }_{1}{  \sigma }_{2}^{-1}$ & $-1$ &
${  1.61803 }$
& ${\{1+t}^{2},...\}$ & ${  5}_{1}^{2}{  \stackrel{5}\sim H}$ &
 \\ \hline
&&&&&&\\
${  (19)=({\overline{18}})}$ & ${  \sigma }_{1}^{-1}{  \sigma }_{2}^{2}{  \sigma }_{1}^{-1}{  \sigma }_{2}$ &
$-1$ & 1.61803  & $\{1+t^2,... \}$ &
${\overline{5_1^2}} \stackrel{5}\sim {  H}
$ &  \\ \hline
  &&&&&&\\
${  (20)=({\overline{16}})}$ & ${  \sigma }_{2}^{-2}{  \sigma }_{1}^{-1}{  \sigma }_{2}{\sigma}_{1}^{-1}$ &
${  1}$ & $%
{  1}$ & ${  \{1,...\}}$ & $4_1^2 \stackrel{5}\sim {  T}_{1}$ &  \\ \hline
&&&&&&\\
${  (21)}$ & ${\sigma}_{2}^{-2}{  \sigma }_{1}^{2}{\sigma }_{2}^{-1}{\sigma}_{1}$ & 1 & ${1.61803 }$ & ${  %
\{1+t^2, ...\}}$ & ${6}_{3} \stackrel{5}{\sim} H$ &  \\ \hline
&&&&&&\\
&&&&&$6_{3}^{3}$ e.g. &\\
${  (22)=({\overline{27}})}$ & ${  \sigma }_{1}^{2}{  \sigma }_{2}{  \sigma }_{1}^{2}{  \sigma }_{2}$ & 1 &
${  2.14896 }$ & $  \{1+t-t^3,... \}$ &  $(3,3)$-torus link &$ ({\sigma }_{1}{  \sigma }_{2}{  )}^{3}$\\
&&&&&or $P_{[2,2,-2]}$ &  \\ \hline&&&&&&\\
&&&&&${\overline{6_3^3}}$ e.g.&\\
${(23)=({\overline{26}})}$ & ${  \sigma }_{1}{  \sigma }_{2}^{-2}{  \sigma }_{1}{  \sigma }_{2}^{2}$ & 1 &
${  2.14896 }$ & $  \{1+t+t^3,... \}$ & (3,-3)-torus link &\\
&&&&&or $P_{[2,-2,-2]}$ &  \\ \hline
&&&&&&\\
&&&&&&${\sigma}_{1}{\sigma}_{2}^{2}{\sigma}_{1}^{2}{  \sigma }_{2}^{-2}={\overline{7_8^2}}$\\
&&&&& ${\overline{6_1^3}}$ rep. $P[-2,-2,-2]$ & ${({\sigma}_{1}{\sigma}_{2})}^{4}$ \\
${(24)=({\overline{25}})}$ & $({\sigma}_{1}{\sigma}_{2}^{-2})^{2}$ & $-1$ & ${  2.49721 }$
& ${\{2+t}^{2}{,... \}}$ &${\overline{6_1^3}} \stackrel{5}\approx
 {\overline{7_8^2}}$ $\stackrel{5}\approx {8_{19}}$   &
 $= \sigma_2\sigma_1 \sigma_2^2\sigma_1^2\sigma_2\sigma_1= 8_{19}$
  \\ \hline
\end{tabular}
\newpage
\begin{center}
\begin{tabular}{|l|l|l|l|l|l|l|}
\hline
${GAP~CC\#}$ & {braids(shortest)} & ${\tiny{{\bf{F}}}}$& ${V(L)}$ & ${V}_{L}{(t,5)}$ {rep}${.}$ & {Link} in
$B_3/(a^5)$ & {interesting rep.} \\ \hline &&&&&&\\
&&&&&   &${  \sigma }_{1}^{-1}{\sigma}_{2}^{-2}{\sigma}_{1}^{-2}{  \sigma }_{2}^{2}={7_8^2}$\\
&&&&&$6_1^3$ rep. $P[2,2,2]$& ${({\sigma}_{1}{\sigma}_{2})}^{-4}$\\
${(25)=({\overline{24}})}$ & ${({\sigma}_{1}^{-1}{\sigma}_{2}^2)}^{2}$ & $-1$ & ${2.49721}$
& ${\{1+2t^2,...\}}$ & ${6_1^3} \stackrel{5}\approx{{7_8^2}}$$\stackrel{5}\approx {\overline{{8}_{19}}}$ &
$= {\sigma}_2^{-1}{\sigma}_1^{-1} {\sigma}_2^{-2}{\sigma}_1^{-2}{\sigma}_2^{-1}{\sigma}_1^{-1}$\\
&&&&&&$={\overline{8_{19}}}$ \\ \hline
&&&&&&\\
&&&&&$6_3^3$ e.g. &\\
${(26)=({\overline{23}})}$ & ${  \sigma }_{1}^{-1}{  \sigma }_{2}^{2}{  \sigma }_{1}^{-1}{  \sigma }_{2}^{-2}$
 & 1 & ${2.14896 }$ & $\{1+t-t^3,... \}$  & $(3,3)$-torus link  &  \\
 &&&&&or $P_{[2,2,-2]}$&\\ \hline
&&&&&&\\
 &&&&&${\overline{6_3^3}}$&\\
${  (27)=({\overline{22}})}$ & ${  \sigma }_{1}^{-2}{  \sigma }_{2}^{-1}{  \sigma }_{1}^{-2}{  \sigma }_{2}^{-1}$
& 1 & ${  2.14896 }$ & $  %
\{1+t+t^3,... \}$ & $(3,-3)$-torus link  & ${  (\sigma }_{1}{  \sigma }_{2}%
{  )}^{-3}$ \\
&&&&&or $P_{[2,-2,-2]}$&\\\hline &&&&&&\\
${(28)=({\overline{29}})}$ & ${\sigma}_{2}^{2}{\sigma}_{1}^{-1}{\sigma}_{2}^{2}{  \sigma }_{1}^{-2}$ &
$-\sqrt{5}$ & ${  1.90211 }$ & $\{1+t,... \}$ & $
7_5^2 \stackrel{5}\sim {\overline{7_5^2}}\stackrel{5}\sim {\overline{7_1^3}}$  &  \\
&&&&& ${\overline{7_1^3}}=P_{[-2,-2,-2,-1]}$  &\\
\hline &&&&&&\\
${  (29)=({\overline{28}})}$ & ${  \sigma }_{2}^{-2}{  \sigma }_{1}{\sigma}_{2}^{-2}{\sigma}_{1}^2$ & $-\sqrt{5}$
& ${  1.90211 }$ & $\{1+t,... \}$ & $
{\overline{7_5^2}}\stackrel{5}\sim {7_5^2} \stackrel{5}\sim {7_1^3}$  &  \\
&&&&& ${7_1^3}=P_{[2,2,2,1]} $&\\ \hline
&&&&&&\\
${(30)}$ & $({\sigma}_{1}{\sigma}_{2}^{-1})^3$ & 1 & ${  3.23607 }$
& ${  \{2+2t^2,... \}}$ & ${6}_{2}^{3} $ (Borromean rings) &  \\
& &&&&&\\ \hline
&&&&&&\\
$(31)=({\overline{32}})$ & $({\sigma }_{1}{\sigma}_{2}^{-1})^3{\sigma }_{2}^{-1}$ & 1 & ${  2.14896 }$ &
${\{1+t-t^3,... \}}$ & ${6_3^3}$  & ${  (\sigma }_{1}{  \sigma }_{2}{  )}%
^{7}$ \\
&&&&&&\\ \hline
&&&&&&\\
${(32)=({\overline{31}})}$ & $({\sigma }_{1}^{-1}{\sigma}_{2})^3{\sigma }_{2}$ & 1 & ${2.14896 }$ &
$\{1+t+t^3,... \}$ & ${\overline{6_3^3}}$ & ${(\sigma }_{1}{  \sigma }%
_{2}{  )}^{-7}$ \\
&&&&&&\\
 \hline
&&&&& &\\
${(33)=({\overline{35}})}$ & $({\sigma}_{2}^{2}{\sigma}_{1}^{2})^2$ & $\sqrt{5}$ & ${1.17557 }$ &
$\{1-t^2,... \}$ & ${\overline{8_{10}^3}} \stackrel{5}\approx {8}_{16} \stackrel{5}\sim {\overline{9_{57}^2}} $
& $({\sigma}_{1}{\sigma}_{2}^{-2})^2{\sigma}_{1}{\sigma}_{2}^{-1}={8}_{16}$ \\
\hline &&&&&&\\
${(34)}$ & $({\sigma}_{1}^{-2}{\sigma}_{2}^2)^2$ & $-1$ & ${0.61803 }$ &
$\{1-t,... \}$ & $8_4^3$ & $({\sigma}_{1}{\sigma}_{2})^{10}$, $({\sigma}_{1}{\sigma}_{2})^{-10}$\\ \hline
&&&&& &\\
${(35)=({\overline{33}})}$ & $({\sigma}_{2}^{-2}{\sigma}_{1}^{-2})^2$ &$\sqrt{5}$  & ${1.17557 }$ &
$\{1-t^2,... \}$  & $8_{10}^3 \stackrel{5}\approx {\overline{8_{16}}} \stackrel{5}\sim 9_{57}^2 $ &
$({\sigma}_{1}^{-1}{\sigma}_{2}^{2})^2{\sigma}_{1}^{-1}{\sigma}_{2}={\overline{8_{16}}}$ \\
\hline &&&&&&\\
&&&&&${{9_{40}^2}}\stackrel{5}\approx {\overline{9_{61}^2}}$&
${\sigma}_2{\sigma}_1^{-1}{\sigma}_2{\sigma}_1^2{\sigma}_2^2{\sigma}_1^2={\overline{9_{61}^2}}$\\
${(36)=({\overline{40}})}$ & $({\sigma}_{1}{\sigma}_{2}^{-2})^3$ & 5
& ${  1.17557 }$ & $\{1-t^2,... \}$ & $\stackrel{(2,2)}\sim 9_{49}$ & ${(\sigma}_{1}{\sigma}_{2})^{6}$ \\
\hline
&&&&&&\\
${(37)={\overline{(38)}}}$ & ${\sigma}_{1}({\sigma}_1{\sigma}_2^{-1})^4$ & 1 & ${1.54335 }$ & $\{1+t-t^2,... \}$
& ${\overline{9_{42}^2}}
\stackrel{5}\approx {9_{41}^2}$&
${\sigma}_1^{-1}{\sigma}_{2}{\sigma}_{1}^{-1}{\sigma}_{2}^{2}{\sigma}_{1}^{-2}{\sigma}_{2}^{2}={9_{41}^2}$ \\
\hline &&&&&&\\
${(38)=({\overline{37}})}$ & ${\sigma}_{1}^{-1}({\sigma}_{1}^{-1}{\sigma}_{2})^4$ &1 & ${1.54335 }$ &
$\{1-t-t^2,... \}$  &
$9_{42}^2 \stackrel{5}\approx {\overline{9_{41}^2}}$  &
${\sigma}_{1}{\sigma}_2^{-1}{\sigma}_{1}{\sigma}_2^{-2}{\sigma}_{1}^2{\sigma}_{2}^{-2}={\overline{9_{41}^2}}$ \\
\hline
&&&&&&\\
${  (39)}$ & $({\sigma}_{1}{\sigma}_{2}^{-1})^4$ & $\sqrt{5}$ & ${2.23607 }$ & ${  \{1+t-t^2-t^3,... \}}$ &
${8}_{18}$ &  \\ \hline
&&&&&&\\
&&&&&${\overline{{9_{40}^2}}} \stackrel{5}\approx {9_{61}^2}$&
${\sigma}_2^{-1}{\sigma}_1{\sigma}_2^{-1}{\sigma}_1^{-2}{\sigma}_2^{-2}{\sigma}_1^{-2}$\\
${(40)=({\overline{36}})}$ & $({\sigma}_{1}^{-1}{\sigma}_{2}^2)^3$ & 5 & ${  1.17557 }$ & $\{1-t^2,... \}$ &
$\stackrel{(2,2)}\sim {9_{49}}$ & $={9_{61}^2}$ \\
&&&&&&${({\sigma}_{1}{\sigma}_{2})}^{-6}$ \\ \hline
&&&&&&\\
&&&&&
&${\sigma}_{1}^{-2}{\sigma}_{2}{\sigma}_{1}^{-2}{\sigma}_{2}^{2}{\sigma}_{1}^{2}{\sigma}_{2}^{2}$\\
&&&&& &$=L11n_{170}$\\
${(41)=({\overline{42}})}$ &  ${\sigma}_{1}{\sigma}_{2}^{-1}{\sigma}
_{1}{\sigma}_{2}^{-2}{\sigma}_{1}{\sigma}_{2}^{-1}{\sigma}_{1}{\sigma}_{2}^{-2}  $    & 1 & ${3.44298 }$ &
$\{1-t-2t^2-t^3,... \}$ &
$L10a163$& ${(\sigma}_{1}{\sigma}_{2}{)}^{9}$ \\
\hline
&&&&&&\\
&&&&&
&${\sigma}_{1}^{2}{\sigma}_{2}^{-1}{\sigma}_{1}^{2}{\sigma}_{2}^{-2}{\sigma}_{1}^{-2}{\sigma}_{2}^{-2}$\\
&&&&& &$={\overline{L11n_{170}}}$\\
${(42)=({\overline{41}})}$ & ${\sigma}_{1}^{-1}{\sigma}_{2}{\sigma}_{1}^{-1}
{\sigma}_{2}^{2}{\sigma}_{1}^{-1}{\sigma}_{2}{\sigma}_1^{-1}{\sigma}_2^2$ &1  & ${  3.44298 }$ &
$\{1+2t+t^2-t^3,... \}$ & ${\overline{L10a163}}$ & ${(\sigma}_{1}{  \sigma }_{2}{)}^{-9}$ \\
\hline
\end{tabular}
\end{center}
\newpage
\begin{tabular}{|l|l|l|l|l|l|l|}
\hline
${GAP~CC\#}$ & {braids(shortest)} & ${\tiny{{\bf{F}}}}$& ${V(L)}$ & ${V}_{L}{(t,5)}$ {rep}${.}$ & {Link} in
$B_3/(a^5)$ & {interesting rep.} \\ \hline
&&&&&&\\
&&&&&&${\sigma}_{1}^{-2}{\sigma}_{2}^{2}{\sigma}_{1}^{-2}{\sigma}_{2}^{-2}{  \sigma }_{1}^{2}{\sigma}_{2}^{-2}$\\
&&&&& ${\overline{11a_{177}}}$&$=\overline{12n_{1958}}$\\
${(43)=({\overline{44}})}$ &
${{\sigma}_1}^{-2}{{\sigma}_2}^{2}{\sigma}_1^{-1}{{\sigma}_2}{\sigma}_1^{-1}{{\sigma}_2}{\sigma}_1^{-1}{\sigma}
_2^{2}$ &5 & ${2.93565 }$ &$\{1-2t^2-t^3,... \}$  &$\stackrel{(2,2)}\sim {\overline{9_{40}}}$   &
${(\sigma}_{1}{\sigma}_{2}{)}^{-12}$ \\
\hline
&&&&&&\\
&&&&&&${\sigma}_{1}^{2}{\sigma}_{2}^{-2}{\sigma}_{1}^{2}{\sigma}_{2}^{2}{\sigma}_{1}^{-2}{\sigma}_{2}^{2}$\\
&&&&&$11a_{177}$&$=12n_{1958} $\\
${(44)=({\overline{43}})}$ &
${{\sigma}_1}^2{{\sigma}_2}^{-2}{\sigma}_1{{\sigma}_2}^{-1}{\sigma}_1{{\sigma}_2}^{-1}{\sigma}_1{\sigma}_2^{-2}$
& 5 & ${2.93565}$ &$\{1+2t-t^3,... \}$  & $\stackrel{(2,2)}\sim {9_{40}}$ & ${(\sigma}_{1}{\sigma}_{2}{)}^{12}$
\\
\hline &&&&&&\\
${(45)}$ & $({\sigma}_{1}{\sigma}_{2}^{-1})^5$& 1 & ${0.381966}$ & $\{1-2t+t^2,...\}$ & ${10}_{123}$ &
${(\sigma}_{2}^2{\sigma}_1^{-2})^{3}={(\sigma}_{1}{\sigma}_{2})^{15}$ \\ \hline
\end{tabular}
}}

\section{$5$-move equivalence of pretzel and Montesinos links}\label{Section 5}

In this section, we deal with classification of pretzel and
Montesinos links up to 5-move equivalence.
The classification is complete for pretzel links and for Montesinos links it
is complete up to an elementary question
(Problem 5.3), having mutation in background\footnote{We do
not deal in this paper with surgery interpretation
of our result, it is worth however to mention that our work can be
related to classifying Seifert fibered manifolds with basis $S^2$
modulo $\pm \frac{1}{5}$-surgeries \cite{D-P-2,D-P-3}.}.
After establishing notation, we formulate the main result of the section.
The proof of Theorem 5.1 is divided into three parts.
First, we identify pretzel and Montesinos link representatives
in 5-move equivalence classes. In Subsection 5.2 we classify
pretzel representatives. In Subsection 5.1 we deal with
Montesinos representatives which are not pretzel links.

Our notation for pretzel and Montesinos links is fairly standard.
It is convenient for us to draw Montesinos links horizontally, so
they look like  pretzel links with columns decorated by rational tangles.
In a pretzel link a column $[n]$ contains $n$ right-handed vertical
half twists (Figures 5.1), in a Montesinos link
$M_{[\frac{p_1}{q_1},...,\frac{p_k}{q_k}]}$ the $i$th column is decorated
by $[\frac{p_i}{q_i}]$ rational tangle (Figure 5.2).
With this notation
we have $P_{[n_1,...,n_k]} = M_{[\frac{1}{n_1},...,\frac{1}{n_k}]}$.
If one column, say $[\frac{p_i}{q_i}]$ is repeated $m$ times in a row,
we write succinctly $M_{[...,m[\frac{p_i}{q_i}],...]}$.
\\ \ \\
\centerline{\psfig{figure=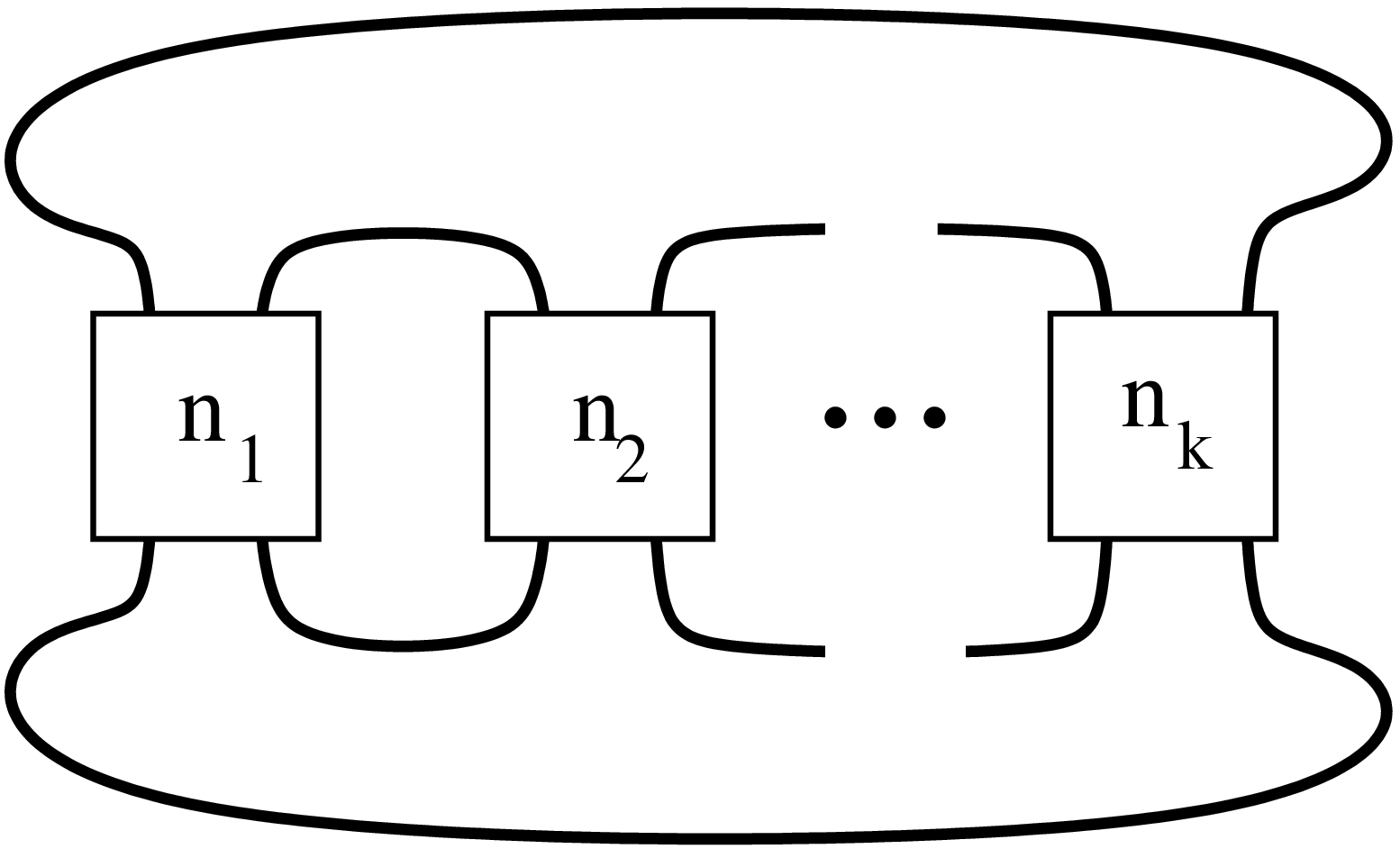,height=2.7cm} \ \ \ \ \ \ \ \
\ \ \ \ \
\psfig{figure=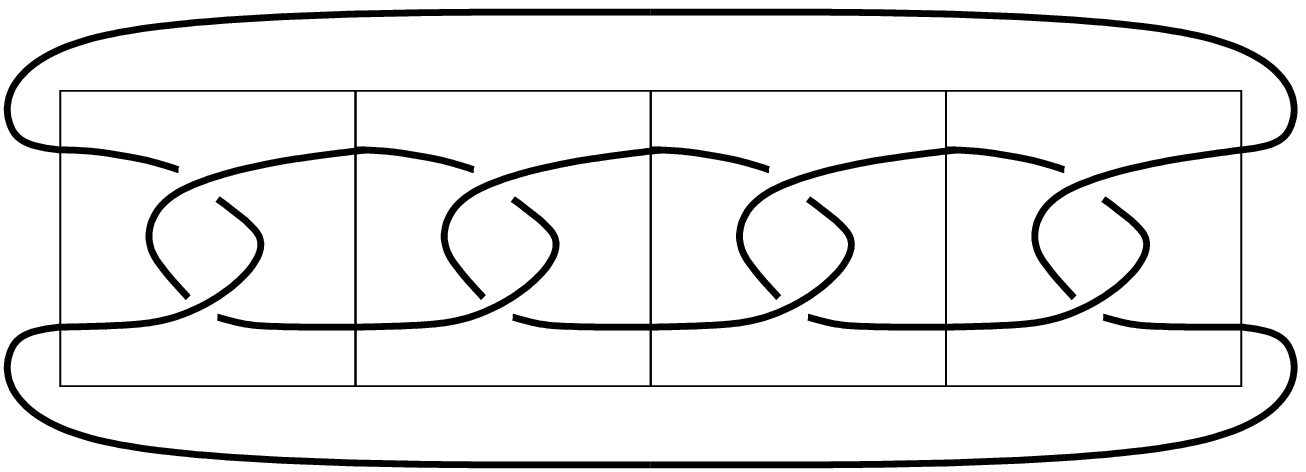,height=2.5cm}}
\begin{center}
Fig. 5.1; $P_{[n_1,...,n_k]}$ \ and \ \ $P_{[2,2,2,2]}$
\end{center}
\ \\ \ \\
\centerline{\psfig{figure=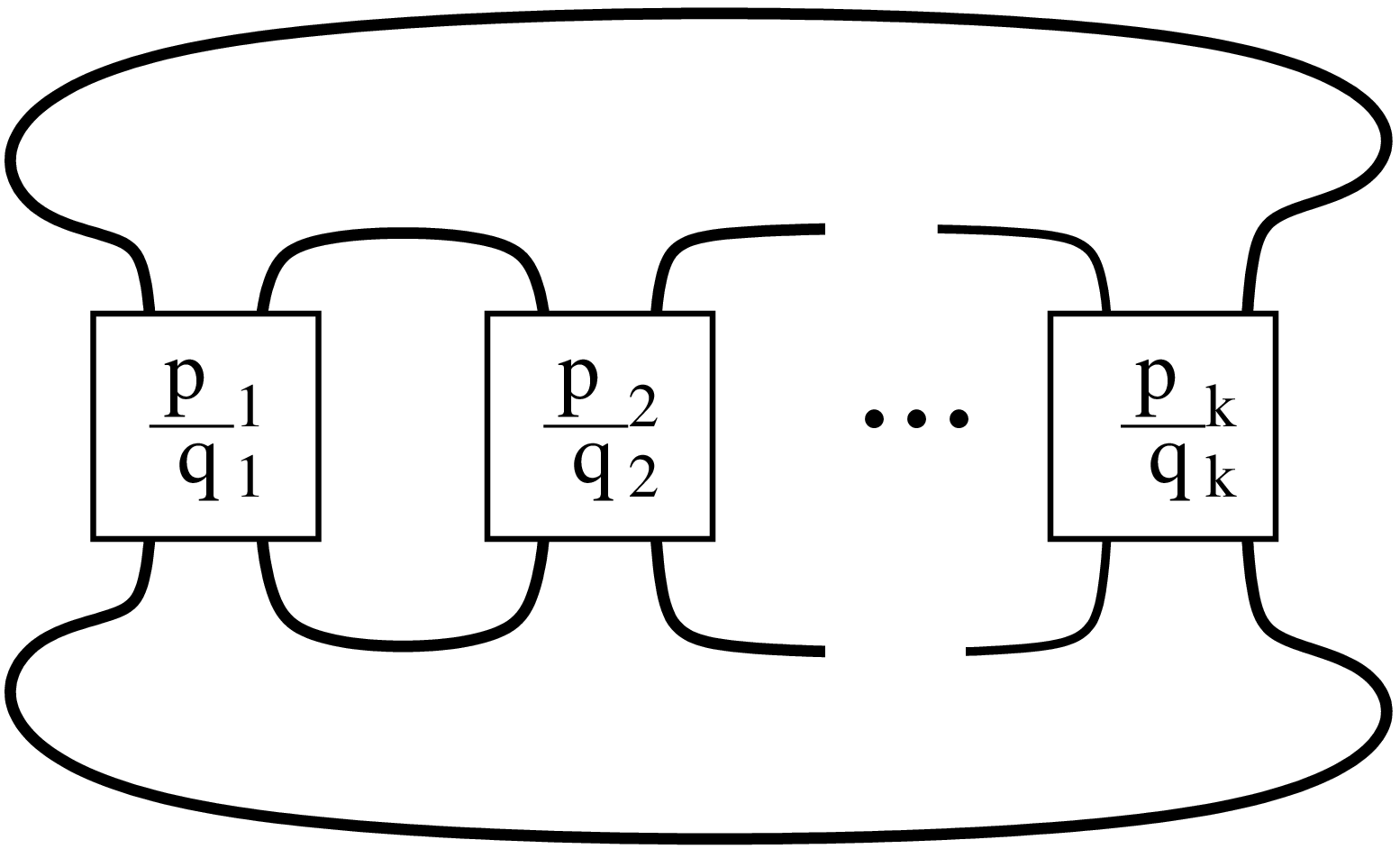,height=2.5cm}
\ \ \ \ \ \ \ \
\psfig{figure=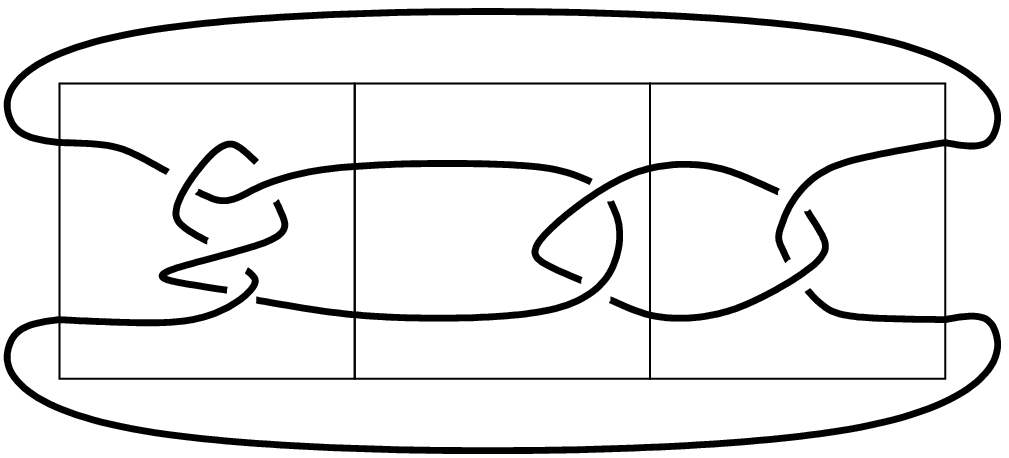,height=2.5cm}}

\begin{center}
Fig. 5.2; $M_{[\frac{p_1}{q_1},...,\frac{p_k}{q_k}]}$\ and \ \
$M_{[\frac{2}{5},\frac{1}{2},\frac{1}{2}]}$
\end{center}

\begin{theorem}\label{Theorem 5.1}\ \\
(1) Every Montesinos link $M_{[\frac{p_1}{q_1},...,\frac{p_k}{q_k}]}$
is $5$-move equivalent to a link from (i), (ii)
or (iii) listed below: \\
(i) Pretzel link $M_{[m[\frac{1}{2}],[s]]}$, for $m\geq 3$,\
 $-2 \leq s \leq 2$.\\
(ii) Montesinos link with all
$\frac{p_i}{q_i}= \frac{2}{5}$ or $\frac{1}{2}$, that is
up to permutation of columns
$M_{[k[\frac{2}{5}],m[\frac{1}{2}]]}$, $k\geq 1$, $k+m \geq 3$.\\
(iii) Connected sum of any number of $T_2$'s, $H$'s or $4_1$'s
(including $T_1$).\\
(2) Links of $M_{[k[\frac{2}{5}],m_1[\frac{1}{2}]]}$, $k\geq 1$, $k+m_1 \geq 3$
and $M_{[m_2[\frac{1}{2}],[s]]}$, for $m_2\geq 3$,\
 $-2 \leq s \leq 2$
are pairwise non $5$-move equivalent and they are not 5-move equivalent
to links listed in (iii) (compare Problem 5.3).
\end{theorem}
\begin{proof}
We prove here part (1) of the theorem. Part (2) will be dealt
with in Subsections 5.2 and 5.3.
Recall that every rational tangle is $5$-move equivalent to one of
the twelve tangles of Lemma 2.10. This is the starting point to
$5$-move classification of Montesinos links. If every column
 $[\frac{p_i}{q_i}]$ of a Montesinos link $M$ is $5$-move equivalent to
a tangle different from $[\frac{2}{5}]$ and $[\frac{1}{0}]$ then $M$ is
$5$-move equivalent to a pretzel link with columns $[\pm\frac{1}{2}]$
or $[\pm 1]$. Furthermore, a column
$[-\frac{1}{2}]$ is isotopic to $[\frac{1}{2}]*[-1]$ and
$[\pm 1]$'s can be collected together, to obtain $M_{[m[\frac{1}{2}],[s]]}$.
Finally $s$ can be reduced modulo 5 by $5$-moves. Notice that for $m\leq 3$
we obtain rational links, as desribed in Example 5.14 (i) and (ii).
We devote Subsection 5.2 to $5$-move classification of pretzel links
$M_{[m[\frac{1}{2}],[s]]}$.

Assume now that at least one column, $[\frac{p_i}{q_i}]$ of $M$
reduces to $[\frac{2}{5}]$ tangle but none to $[\frac{1}{0}]$
(compare Proposition 2.11).
As we checked already when classifying rational tangles up to 5-moves,
$[\frac{2}{5}] \stackrel{5}\sim [\frac{2}{5}\pm1]$ and
$[\frac{2}{5}]*[-\frac{1}{2}] = [\frac{2}{5}]*[\frac{1}{2}]*[-1]
\stackrel{5}\sim [\frac{2}{5}]*[\frac{1}{2}]$, therefore $M$
reduces to a Montesinos link with all
$\frac{p_i}{q_i}= \frac{2}{5}$ or $\frac{1}{2}$ which in fact, after
permutation of columns gives $M_{[k[\frac{2}{5}],m[\frac{1}{2}]]}$, $k\geq 1$.
Notice that $M_{[\frac{2}{5}]} = H$,  $M_{[\frac{2}{5},\frac{1}{2}]}=
[\frac{9}{4}]^N \stackrel{5}\sim T_1$, $M_{[2[\frac{2}{5}]]} =
[\frac{20}{9}]^N \stackrel{5}\sim  T_2$. We devote Subsection 5.1 to
$5$-move classification of Montesinos links
$M_{[k[\frac{2}{5}],m[\frac{1}{2}]]}$.

Finally, if there is a column, say $[\frac{p_i}{q_i}]$, which reduces
to $[\frac{1}{0}]$ tangle then
$M_{[\frac{p_1}{q_1},...,\frac{p_i}{q_i},...,\frac{p_k}{q_k}]}=
[\frac{p_1}{q_1}]^D \# \cdots \# [\frac{p_{i-1}}{q_{i-1}}]^D \#
[\frac{p_{i+1}}{q_{i+1}}]^D\# \cdots \# [\frac{p_k}{q_k}]^D$ and any link
$[\frac{p}{q}]^D$ is $5$-move equivalent to $T_1$, $T_2$, $H$, or $4_1$.
  The proof of Theorem 5.1(1) is completed.
\end{proof}
\begin{remark}\label{Remark 5.2}
The transposition of neighboring columns in a Montesinos link is a mutation
 and cannot be detected by invariants we introduced. The smallest examples
of Montesinos links for which we do not know whether they are $5$-move
equivalent are $12$ crossing, $2$-component links
$M_{[2[\frac{2}{5}],2[\frac{1}{2}]]}$ and
$M_{[\frac{2}{5},\frac{1}{2},\frac{2}{5},\frac{1}{2}]}$, see
Figure 5.3.
\end{remark}
\ \\
\centerline{\psfig{figure=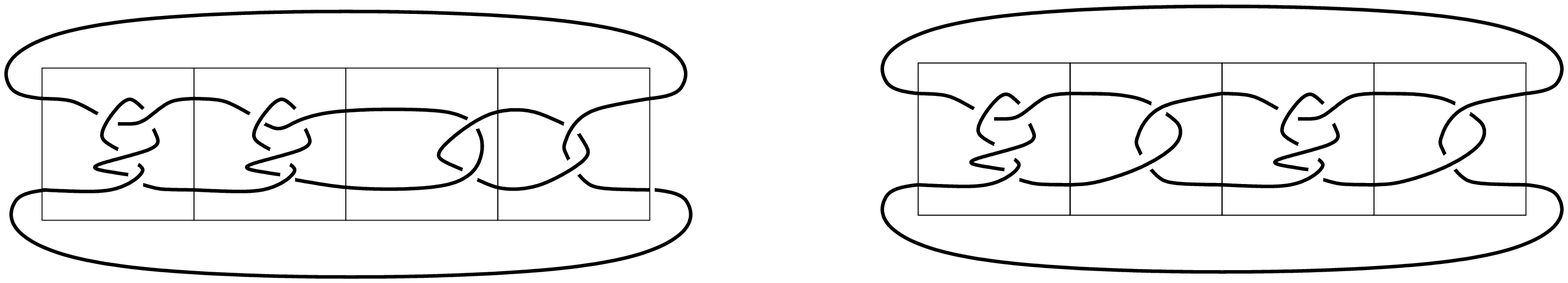,height=2.5cm} }
\ \\
\begin{center}
Fig. 5.3; $M_{[2[\frac{2}{5}],2[\frac{1}{2}]]}$\ \ \ \ and \ \ \ \
$M_{[\frac{2}{5},\frac{1}{2},\frac{2}{5},\frac{1}{2}]}$
\end{center}

More generally we have the following unresolved cases
concerning classification of pretzel and Montesinos links up to
$5$-move equivalence.
\begin{problem}\label{Problem 5.3}
Consider two Montesinos links $L_1$ and $L_2$ both of them with $k \geq 2$
columns $[\frac{2}{5}]$ and $m \geq 2$ columns $[\frac{1}{2}]$
(in any order). Are $L_1$ and $L_2$ $5$-move equivalent?
\end{problem}

The next problem, which we partially solve in Lemma 5.5,
 is related to the possibility that a column of
a Montesinos link is $5$-move equivalent to $[\frac{1}{0}]$ and some other
columns to $[\frac{2}{5}]$ .
\begin{problem}\label{Problem 5.4}\ \\
Let two links $L_1$ and $L_2$ be connected sums of $k_i$ ($k_i\geq 1$)
 copies of $4_1$, $m_i$ copies of $H$, and $n_i$ copies of $T_2$
(taken in any manner).
Are $L_1$ and $L_2$ $5$-move equivalent?\\
Notice that $L_i$ is $(2,2)$-move equivalent to $T_{k_i+n_i+1}$ by
$k_i+m_i$ $(2,2)$-moves thus we can limit the problem to the case
when $k_1+n_1+1 = k_2+n_2+1$ and $k_1 + m_1 \equiv k_2 + m_2 \mod 2$,
compare Remark 3.8, and the last paragraph of Section 3.
\end{problem}

\begin{lemma}\label{Lemma 5.5}
(i)
Let two links $L_1$ and $L_2$ be connected sums of $m$ copies of $H$ and
$n$ copies of $T_2$ (taken in any manner).
Then $L_1$ and $L_2$ are $5$-move equivalent.\\
(ii) Let two links $L_1$ and $L_2$ be connected sums of $k$ copies of $4_1$,
$m$ copies of $H$, and $n$ copies of $T_2$ (taken in any manner).
Then $L_1$ and $L_2$ are $5$-move equivalent.
\end{lemma}
\begin{proof} The main idea is that a 5-move allows us to change the
disjoint sum into connected sum.
We illustrate Lemma 5.5 by an example: the link $H\#H \sqcup T_1$ is
$5$-move equivalent to the disjoin sum $H \sqcup H$.
Namely, by one 5-move we can change $H\#H \sqcup T_1$ to $H\#H\#5_1$.
Similarly, $H \sqcup H$ can be changed by one 5-move to $H\#5_1\#H$.
Since we can choose the connected sum formation in such a way that
$H\#H\#5_1$ and $H\#5_1\#H$ are ambient isotopic (see Figure 5.4),
hence Lemma 5.5 follows in this case. In the case of two different
formations of a connected sum $H\#H\#H$, the 5-move equivalence is
illustrated in Figure 5.5.
The general proof follows the same idea. Similarly
one proves part (ii) of the lemma.
\end{proof}

\ \\
\centerline{\psfig{figure=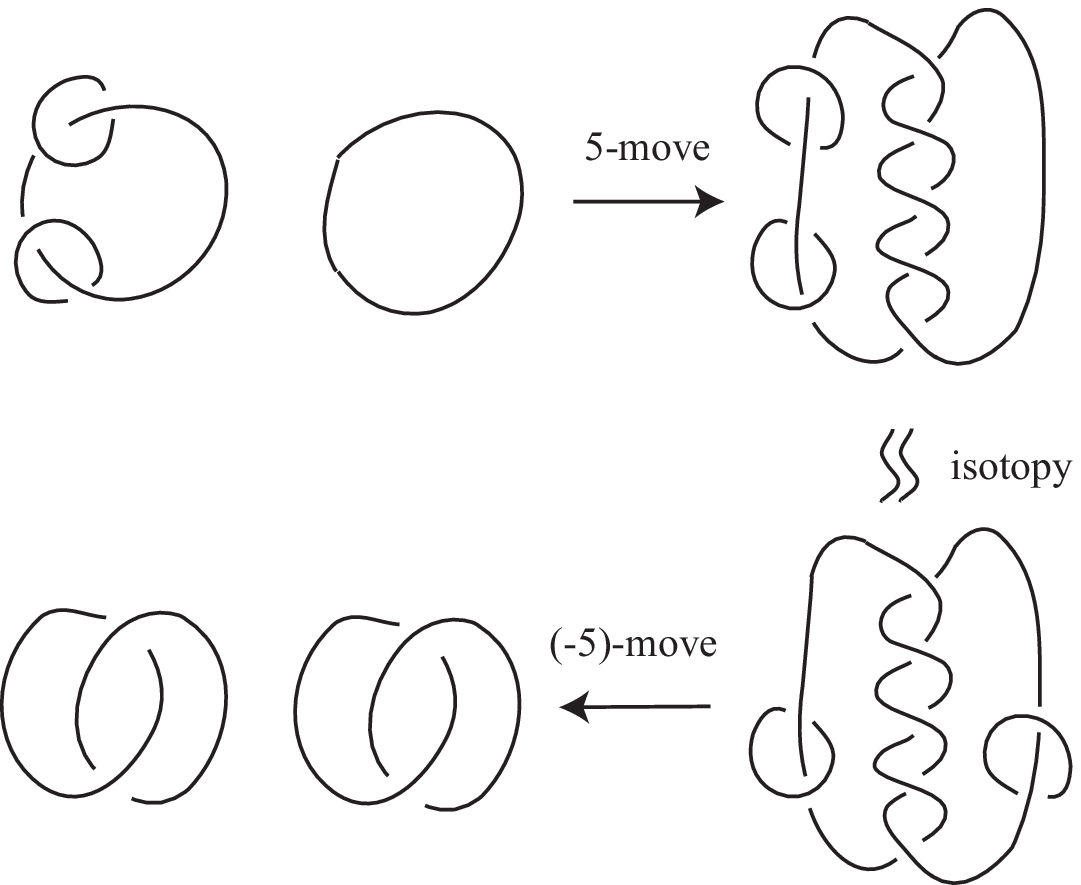,height=6cm}}
\begin{center} Fig. 5.4; $H\#H \sqcup T_1$ and $H \sqcup H$ are
related by two $(\pm 5)$-moves \end{center}

\ \\
\centerline{\psfig{figure=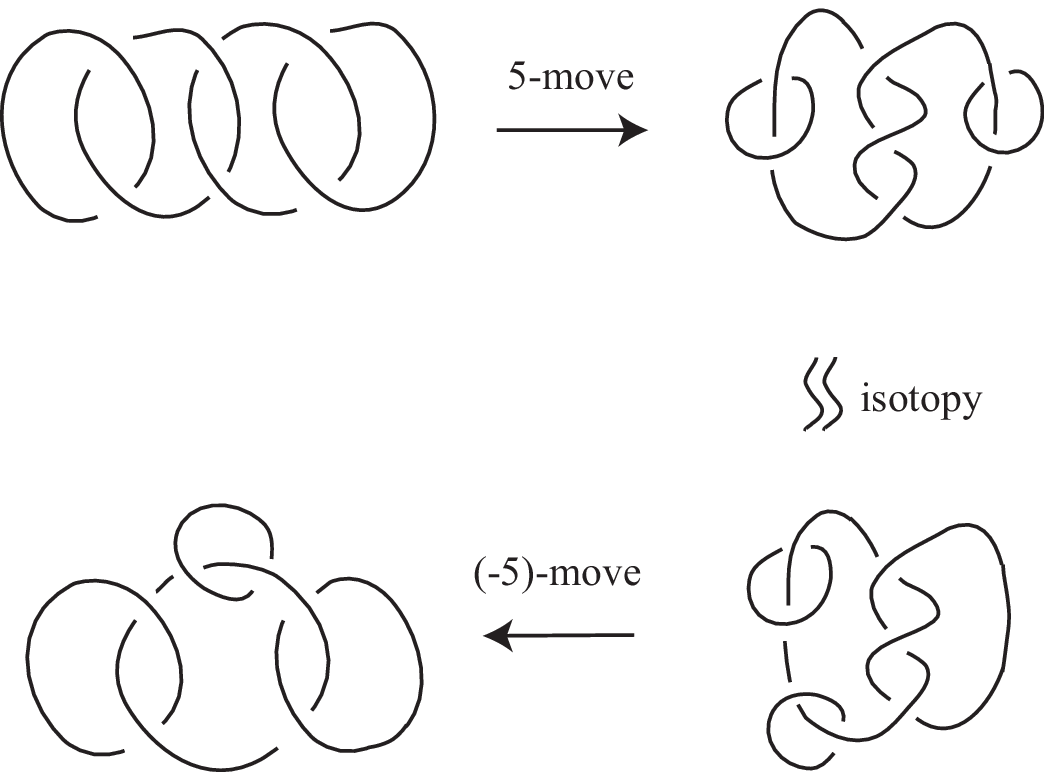,height=5.5cm}}
\begin{center} Fig. 5.5; Two different realizations of $H\#H\#H$ are
related by two $(\pm 5)$-moves \end{center}

\par

\subsection{Kauffman bracket 5-move invariant for
Montesinos links}\label{Subsection 5.1}

In this subsection we compute the Jones (and Kauffman bracket)
polynomial of Montesinos links $M_{[k[\frac{2}{5}],m[\frac{1}{2}]]}$,
$k \geq 1$.
We show that $V(M)= |V_M(e^{\pi i/5})|$ is sufficient to  separate
$5$-move equivalence classes of these links. Our computation is helped by
the fact that
$V_{4_1}(t)= t^{-2}\frac{t^5+1}{t+1}$ so $V([\frac{2}{5}]^D)=V(4_1)=0$,
and the fact we already used that
$[\frac{2}{5}+1] \stackrel{5}\sim [\frac{2}{5}]
\stackrel{5}\sim [-\frac{2}{5}]$.

\begin{example}\label{Example 5.5}
The (prime) Montesinos links, which are not pretzel (or rational)
links, with no more than 8 crossings and up to the mirror image are
$8^2_9=M_{[\frac{2}{5},\frac{1}{2},\frac{1}{2}]}$,
$8^2_{10}=M_{[\frac{3}{5},\frac{1}{2},\frac{1}{2}]}$,
$8^2_{15}=M_{[\frac{2}{5},\frac{1}{2},-\frac{1}{2}]}$  and
$8^2_{16}=M_{[\frac{3}{5},\frac{1}{2},-\frac{1}{2}]}$ (Fig.5.6).
All these links are $5$-move equivalent by identities
$[\frac{3}{5}] \stackrel{5}\sim [\frac{2}{5}]$,
$[\frac{2}{5}]*[-\frac{1}{2}] = [\frac{2}{5}]*[-1]*[\frac{1}{2}]
\stackrel{5}\sim [\frac{2}{5}]*[\frac{1}{2}]$. Observe that
$V_{8^2_9}(t) \equiv t^{-1/2 \pm 2}(1-t) \mod \frac{t^5+1}{t+1}$
(compare Theorem 5.7).
\end{example}
\ \\
\centerline{\psfig{figure=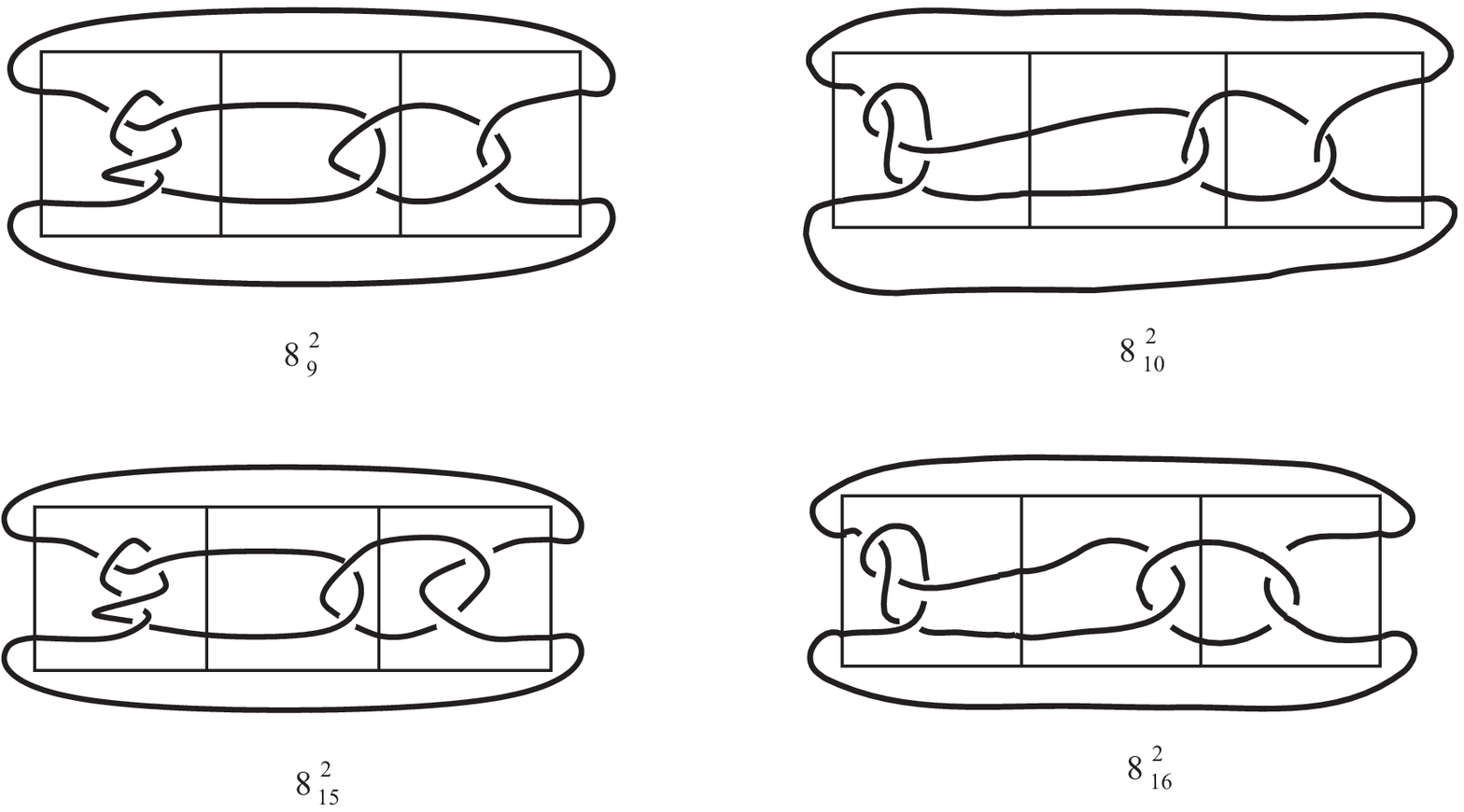,height=8.5cm}}
\begin{center} Fig. 5.6; 5-move equivalent Montesinos links \end{center}

To formulate succinctly the main result of this subsection recall that
$I_t$ denotes the ideal in $Z[t^{\mp \frac{1}{2}}]$ generated by
$\frac{t^5+1}{t+1}$. Let $\doteq$ denote equivalence
up to $\pm t^{i/2}$ for some $i$. Similarly, let $I_A$ be the ideal
in $Z[A^{\mp 1}] $ generated by $\frac{A^{20}+1}{A^4+1}$.

Then for the Jones polynomial modulo $I_t$ we obtain the following
theorem which is the main tool to classify Montesinos links
$M_{[k[{\frac{2}{5}}],m[\frac{1}{2}]]}$ for $k \geq 1$,
up to $5$-move equivalence.


\begin{theorem}\label{Theorem 5.7}
(i)
$V_{M_{[k[\frac{2}{5}],m[\frac{1}{2}]]}}(t)  \et $
$(1+t^2)(1-t^2)^{k-1} (1-t)^{m}$ for
$m \geq 0,k \geq 1$.\\
(ii)If additionally, $k+m \geq 2$ we can write succinctly:
\par
$V_{M_{[k[\frac{2}{5}],m[\frac{1}{2}]]}}(t) \et
(1+t)^{k-1}(1-t)^{k+m-2}$.
\end{theorem}
\ \\ \ \ \ \
{\it Proof:\ \ }
The main observation leading to the proof is that $V_{4_1}(t) \equiv 0$
mod $I_t$.
As before, let $T_A * T_B=$
(\ \parbox{2.5cm}{\psfig{figure=TangleProdDIP.eps,height=0.8cm}}),
and $T_A^N=$  \ \parbox{1.1cm}{\psfig{figure=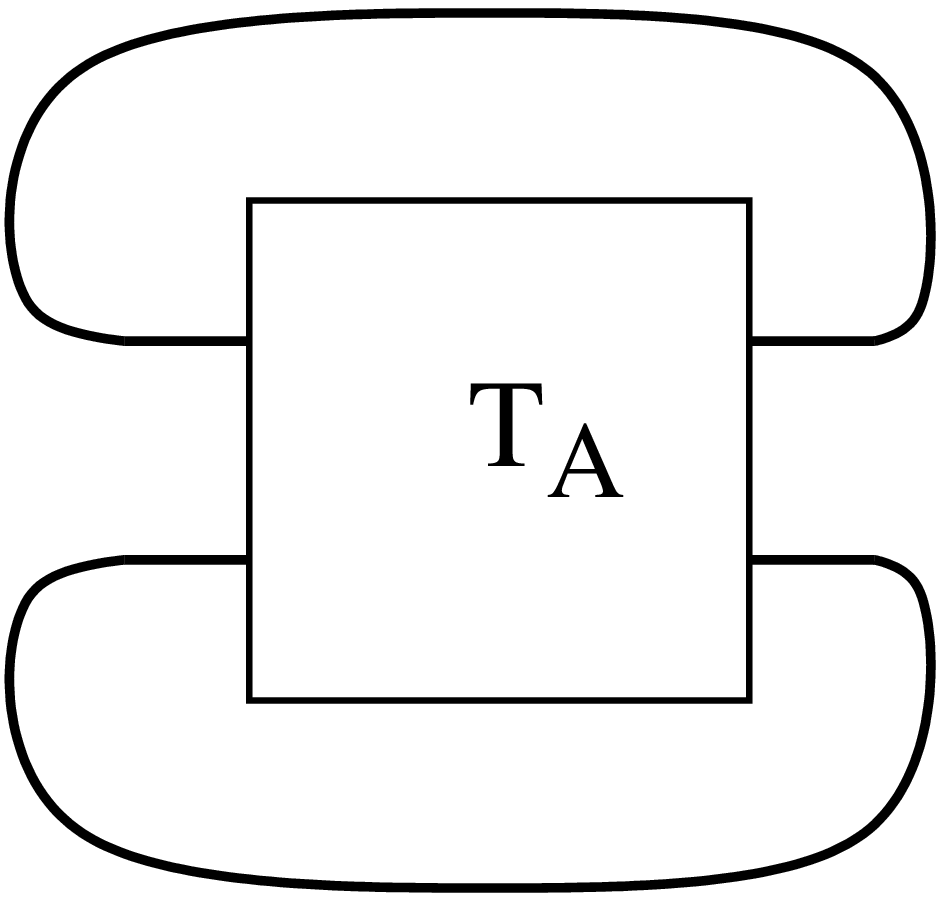,height=0.8cm}}, \
and $T_A^D =$
\ \parbox{1.6cm}{\psfig{figure=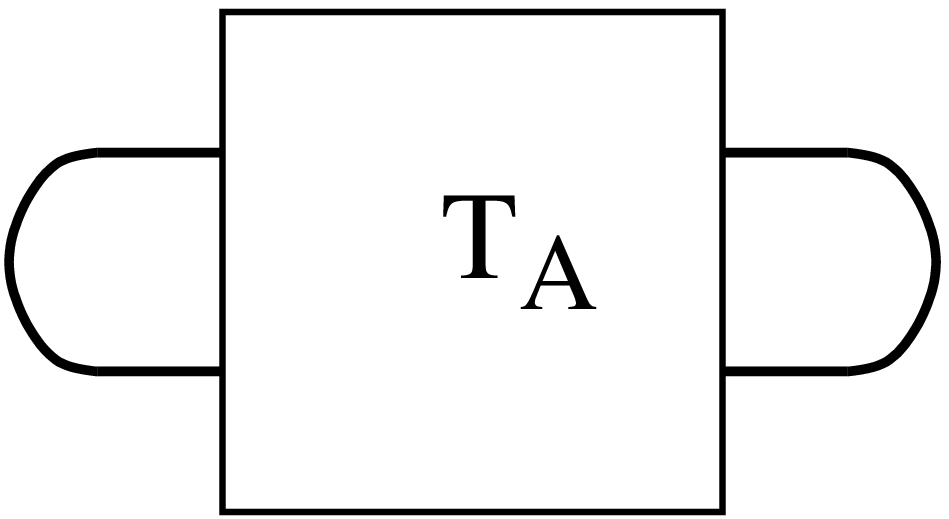,height=0.8cm}}.
 We have the following formulas
for the Kauffman bracket (Lickorish-Millett generalization
of the Conway's formula).


\begin{lemma}\label{Lemma 5.8}
(a) $\langle (T_A*T_B)^D\rangle =\langle  T_A^D\rangle\langle
T_B^D\rangle$. \\
(b) $\langle (T_A*T_B)^N\rangle  =\frac{1}{d^2-1}
(d\langle  T_A^N\rangle\langle  T_B^N\rangle+d\langle  T_A^D\rangle\langle
T_B^D\rangle-\langle  T_A^N\rangle\langle  T_B^D\rangle-\langle
T_A^D\rangle\langle  T_B^N\rangle)$, \\
where $d$ denotes the value of bracket for $T_2$, that is $d=-A^2-A^{-2}$.
\end{lemma}

\par

\vspace{0.5cm}

We use variations of the Conway-Lickorish-Millett formula and we
develop them in the language of the Kauffman bracket skein
modules \cite{Pr-2,H-P}.\\
 The tangles $T_A$ and $T_B$ can be written in a basis of a 2-tangle,
$e_h=\asymp $
and $e_v=)( $ as $T_A=a_1e_h+a_2e_v$. Then $ \langle  T_A^N\rangle=da_1+a_2$,
and $\langle  T_A^D\rangle=a_1+da_2$.
Similarly $T_B=b_1e_h+b_2e_v, \langle  T_B^N\rangle=db_1+b_2$,
and $\langle  T_B^D\rangle=b_1+db_2$, $a_1,a_2,b_1,b_2 \in Z[A^{\mp}]$.\\
From this we have:\\
$(d^2-1)a_1= d\langle  T_A^N\rangle-\langle  T_A^D\rangle$,\ \
$(d^2-1)a_2= d\langle  T_A^D\rangle-\langle  T_A^N\rangle$,\\
$(d^2-1)b_1= d\langle  T_B^N\rangle-\langle  T_B^D\rangle$, and
$(d^2-1)b_2= d\langle  T_B^D\rangle-\langle  T_B^N\rangle$.

Finally we get:
\begin{lemma}\label{Lemma 5.9}
(i) $T_A*T_B=a_1b_1e_h+(a_1b_2+a_2b_1+a_2b_2d)e_v$, \\
(ii) $\langle (T_A*T_B)^D\rangle=a_1b_1+(a_1b_2+a_2b_1+a_2b_2d)d=
(a_1+a_2d)(b_1+b_2d)$, and \\
(iii) $\langle (T_A*T_B)^N\rangle = (a_1b_1+a_2b_2)d+a_1b_2+a_2b_1 =
a_1(b_1d+b_2)+a_2(b_1+b_2d) = \\
 a_1\langle  T_B^N\rangle+a_2\langle  T_B^D\rangle$.
\end{lemma}

Formula (iii) leads immediately to the formula of Lemma 5.7(b).

\begin{example}\label{Example 5.10}
 In the Kauffman bracket skein module of $2$-tangles
we get the following: \
$\langle$ $\parbox{0.6cm}{\psfig{figure=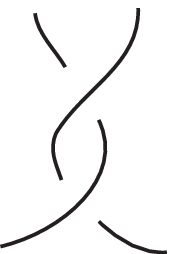,height=0.8cm}}
\rangle = (1-A^{-4}) \langle \sg \rangle $
$+ A^2 \langle )( \rangle $; \ that is
$[\frac{1}{2}]=(1-A^{-4})e_h +A^2e_v$\\
$\langle \parbox{0.8cm}{\psfig{figure=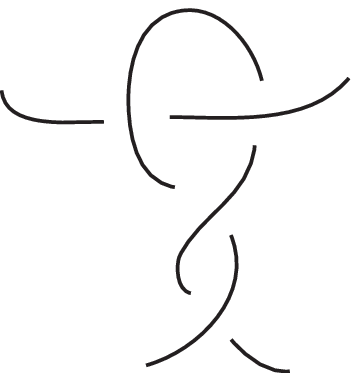,height=0.8cm}}
\rangle $
$=(A^{-8}-A^{-4}+2-A^4) \langle \sg \rangle $
$+(A^2-A^6)  \langle )(\rangle$,\\
that is $[\frac{2}{5}]=(A^{-8}-A^{-4}+2-A^4)e_h + (A^2-A^6)e_v$\\
In particular, \
$a_1([\frac{2}{5}])=A^{-8}-A^{-4}+2-A^4 $
$=_{t=A^{-4}} t^2-t+2-t^{-1}$ $\et 1-t^2$.
\end{example}

Notice that $(t^2-1)(t^2+1)=(t+1)(t-1)(t^2+1)=(t+1)(t^3-t^2+t-1)
\et t^4(1+t)$; in particular, $1+t^2 \et (1-t)^{-1}$ the observation
used to derive (ii) from (i) in Theorem 5.7.
\begin{lemma}\label{Lemma 5.11}
(i)  $\langle (T_A*([\frac{2}{5}])^N\rangle \stackrel{I_A}\equiv
a_1A^{-6}\langle  H\rangle=a_1(-A^{-2})(1+A^{-8})$\\
(ii) For $k\geq 1$, $\langle  M_{[k[\frac{2}{5}],m[{\frac{1}{2}}]]}\rangle $
$\ea (-A^{-2})(1+A^{-8})(A^{-8}-A^{-4}+2-A^4)^{k-1}(1-A^{-4})^m$,\\
(iii) $\langle  M_{[k[\frac{2}{5}],m[{\frac{1}{2}}]]}\rangle_{t=A^{-4}}
\et -t^{\frac{1}{2}}(1+t^2)(t^2-t+2-t^{-1})^{k-1}(1-t)^m \et$\\
\ \ \ \ \ \ \ \ \ \ \ \ \ $(1+t^2)(1-t^2)^{k-1}(1-t)^m$.
\end{lemma}

\begin{proof}
By Lemma 5.9(iii)
$\langle (T_A*([\frac{2}{5}]))^N\rangle
= a_1 \langle ([\frac{2}{5}])^N\rangle +a_2 \langle ([\frac{2}{5}])^D\rangle $\\
$\stackrel{I_A}\equiv  a_1 \langle ([\frac{2}{5}])^N\rangle +a_2 \cdot 0
=a_1A^{-6}\langle  H\rangle =a_1 A^{-6}(-A^4-A^{-4})$, and Lemma 5.11(i)
follows.\\
We can write $M_{[k[\frac{2}{5}],m[{\frac{1}{2}}]]}$ as
 $((k-1)[\frac{2}{5}]*m[{\frac{1}{2}}]*[\frac{2}{5}])^N$ and use the previous
formula for
$T_A=(k-1)[\frac{2}{5}]*m[{\frac{1}{2}}]$,  we have in this case
$a_1=(A^{-8}-A^{-4}+2-A^4)^{k-1}(1-A^{-4})^m$, and Lemma 5.11(ii) follows. \
Lemma 5.11(iii) follows after substituting $t=A^{-4}$ and a
 simple calculation modulo $I_t$.

\end{proof}

\begin{example}\label{Example 5.12}
$V_{M_{[\frac{2}{5},\frac{1}{2},\frac{2}{5},\frac{1}{2}]}}(t) =
V_{M_{[2[\frac{2}{5}],2[\frac{1}{2}]]}}(t) \et (1-t^2)(1-t)= 1-t-t^2+t^3$.
\end{example}

We can use Theorem 5.7 to prove a part of Theorem 5.1(2). That is:
\begin{corollary}\label{Corollary 5.13}\ \\
(i) Montesinos links $M_{[k[\frac{2}{5}],m[\frac{1}{2}]]}$,
$k \geq 1$, $m \geq 0$, are pairwise not $5$-move equivalent.\\
(ii) $M_{[\frac{2}{5}]}= H$,
$M_{[\frac{2}{5},\frac{1}{2}]}=
[\frac{9}{4}]^N \stackrel{5}\sim T_1$, $M_{[2[\frac{2}{5}]]} =
[\frac{20}{9}]^N \stackrel{5}\sim  T_2$.\\
(iii) $M_{[k[\frac{2}{5}],m[\frac{1}{2}]]}$, $k \geq 1$, $m+k \geq 3$,
is not $5$-move equivalent to disjoint or connected sums of $T_1$, $H$ or
$4_1$.
\end{corollary}
\begin{proof}
From Theorem 5.7 and Lemma 6.1 it follows that
$V_{M_{[k[\frac{2}{5}],m[\frac{1}{2}]]}}(t)  \et
V_{M_{[k'[\frac{2}{5}],m'[\frac{1}{2}]]}}(t)$ iff $k=k'$ and $m=m'$.
In fact it follows from the proof of Lemma 6.1 that
 values of $V(M_{[k[\frac{2}{5}],m[\frac{1}{2}]]})$, $k\geq 1$,
 are all different (and different from $0$).\ To prove (iii), first
notice that if a link $L$ has $4_1$ as a connected or disjoint sum
 summand, then $V_L(t) \et 0$. Furthermore, if $L$ is a finite
 disjoint or connected sum of $T_1$ and $H$ then $<L> =
(-A^2-A^{-2})^i(-A^4-A^{-4})^j$ and $V_L(t)\et (1+t)^i(1+t^2)^j$.
Therefore $V_L(t)\et V_{M_{[k[\frac{2}{5}],m[\frac{1}{2}]]}}(t)$ only
for $k=1,m=0$ or $k=m=1$, or $k=2,m=0$, as described in (ii) of
Corollary 5.13.
\end{proof}

\par\vspace{1cm}

\subsection{Jones polynomial for pretzel links}\label{Subsection 5.2}
We develop here formulas for Jones polynomial and
Kauffman bracket sufficient to distinguish 5-move
equivalence classes of pretzel links $P_{[2,...,2,1,...,1]}$ and to
complete the classification of pretzel links up to $5$-move equivalence.
To be consistent with notation for Montesinos links we will write
$M_{[m[\frac{1}{2}],[s]]}$ for the pretzel link with $m$ two's and $s$
one's. We have shown in Theorem 5.1 that any pretzel link is $5$-move
equivalent to one of $M_{[m[\frac{1}{2}],[s]]}$ for $-2 \leq s \leq 2$,
or to the connected sum of $T_2$'s and $H$'s.
We will show in this subsection that for $m \geq 3$ such links are not
$5$-move equivalent.
\begin{example}\label{Example 5.14}\
(i) For $m<3$, links $M_{[m[\frac{1}{2}],[s]]}$ are rational links.
Specifically, $M_{[\frac{1}{2}]} = M_{[[\frac{1}{2}],[-1]]} = T_1$,
$M_{[\frac{1}{2},1]} = \bar 3_1 \stackrel{5}\sim H
\stackrel{5}\sim 3_1 = M_{[[\frac{1}{2}],[-2]]}$,
$M_{[[\frac{1}{2}],[2]]} = 4_1$, and \\
$M_{[2[\frac{1}{2}]]} = \bar 4^2_1 \stackrel{5}\sim  T_1
\stackrel{5}\sim 4^2_1 =  M_{[2[\frac{1}{2}],[-2]]} $,
$M_{[2[\frac{1}{2}],[1]]} \stackrel{5}\sim  H
\stackrel{5}\sim   M_{[2[\frac{1}{2}],[2]]}$,
$M_{[2[\frac{1}{2}],[-1]]} = T_2$.\\
(ii) We list here links $M_{[m[\frac{1}{2}],[s]]}$ for $m=3$,
 $-2 \leq s \leq 2$, with their
5-move invariants sufficient
to separate them among themselves and from rational links.\\
  $M_{[3[\frac{1}{2}]]} = 6^3_1$, with $V(6^3_1)\approx 2.49721$,
$V_{6^3_1}(t,5) = \{2+t^2,...\}$,\\
$M_{[3[\frac{1}{2}],[2]]} \stackrel{5}\sim \bar 6^3_1$,
with $V(\bar 6^3_1)\approx 2.49721$, $V_{\bar 6^3_1}(t,5) = \{1+2t^2,...\}$,\\
$M_{[3[\frac{1}{2}],[1]]} = 7^3_1 \stackrel{5}\sim \bar 7^3_1$,
with $V(7^3_1)\approx 1.90211$, $V_{7^3_1}(t,5) = \{1+t,...\}$,
$F_{7^3_1}(1,2cos 2\pi/5)=-\sqrt{5}$, \\
$M_{[3[\frac{1}{2}],[-1]]} = 6^3_3$,
with $V(6^3_3)\approx 2.14896$, $V_{6^3_3}(t,5) = \{1+t-t^3,...\}$,\\
$M_{[3[\frac{1}{2}],[-2]]} = \bar6^3_3$,
with $V(\bar 6^3_3)\approx 2.14896$, $V_{6^3_3}(t,5) = \{1+t+t^3,...\}$.\\
(iii) For $m=4$, the invariant $V_L(t,5)$ separates links:\\
$M_{[4[\frac{1}{2}]]} = 8^4_1$, with $V(8^4_1)\approx 3.67044$,
$V_{8^4_1}(t,5) = \{1+2t+2t^3,...\}$,\\
$M_{[4[\frac{1}{2}],[1]]} = 9^4_1 \stackrel{5}\sim \bar 8^4_1$,
with $V(\bar 8^4_1) \approx 3.67044$, $V_{\bar 8^4_1}(t,5) =
\{2+2t^2+t^3,...\}$,\\
$M_{[4[\frac{1}{2}],[-1]]} = 8^4_2$,
with $V(8^4_2)\approx 3.44298$, $V_{8^4_2}(t,5) = \{3+t^2,...\}$,\\
$M_{[4[\frac{1}{2}],[2]]} \stackrel{5}\sim \bar 8^4_2$,
with $V(\bar 8^4_2)\approx 3.44298$, $V_{\bar 8^4_2}(t,5) = \{1+3t^2,...\}$,\\
$M_{[4[\frac{1}{2}],[-2]]} = 8^4_3$,
with $V(8^4_3)\approx 3.80423$, $V_{8^4_3}(t,5) = \{2+2t,...\}$.
\end{example}

Our main tool to separate links $M_{[m[\frac{1}{2}],[s]]}$ is the
Jones polynomial (or the Kauffman bracket).
\begin{proposition}\label{Proposition 5.15} \ \\
(i) $\langle M_{[m[\frac{1}{2}]]} \rangle =
(1-A^{-4})^md+\frac{(-A^4-A^{-4})^m-(1-A^{-4})^m}{d}$\\
$\stackrel{t=A^{-4}}{=} -(1-t)^m(t^{1/2} + t^{-1/2}) +
\frac{(-1)^{m-1}(t+t^{-1})^m + (1-t)^m}{t^{1/2} + t^{-1/2}}$.\\
(ii) $\langle M_{[m[\frac{1}{2}],[s]]} \rangle =
(-A^3)^s\bigl(\frac{(-A^4-A^{-4})^m}{d} +
(1-A^{-4})^m\bigl(-A^{-4})^s(d -d^{-1})
 \bigr).$\\
(iii) In other words
$$\tilde V_{M_{(m[\frac{1}{2}],[s])}}(t) =
(-1)^{m-1}\frac{(t+t^{-1})^m}{t^{1/2} + t^{-1/2}}
- (1-t)^m\bigl((-t)^s\frac{t+1+t^{-1}}{t^{1/2} + t^{-1/2}}\bigr).$$

\end{proposition}

\begin{proof}
Let $T$ be any 2-tangle,  and let $T^{(m)}$ denote $T \ast \cdots \ast T$.
In the Kauffman bracket skein module we write:
 $\langle  T \rangle =a_1 \langle  \sg \rangle $
 $+a_2 \langle  )( \rangle$, \ and
 $\langle  T^{(m)} \rangle =a_1^{(m)} \langle  \sg \rangle $
$+a_2^{(m)} \langle  )( \rangle $. As before, we have
 $\langle  T^D\rangle=a_1 +a_2d$,
 $\langle  T^N \rangle =a_1 d +a_2 $, and
$\langle (T^{(m)})^N \rangle
=a_1^{(m)} d + {a_2}^{(m)}$, $\langle (T^{(m)})^D \rangle
=\langle  T^D \rangle^m =a_1^{(m)}  + d {a_2}^{(m)}$. Then \
$\langle T^{(m)} {\rangle} = (a_1 \langle  \sg \rangle
+a_2 \langle  )( \rangle )
\ast (a_1 \langle  \sg \rangle
+a_2 \langle  )( \rangle )
\ast \cdots \ast (a_1 \langle  \sg \rangle
+a_2 \langle  )( \rangle )
=a^{(m)} \langle  \sg \rangle + a_2^{(m)} \langle )( \rangle $.
First, we conclude that $a_1^{(m)}=a_1^m$, then $da_2^{(m)}=
\langle  (T^{(m)})^D\rangle -a_1^m$ and from this
$\langle  (T^{(m)})^N \rangle=a_1^m d
+(\langle  T^D \rangle^m -a_1^m)d^{-1} $. \
Specifically for $T=[\frac{1}{2}]=
\parbox{0.6cm}{ \psfig{figure=rational-12DIP.eps,height=0.4cm}} $,
we have  $\langle T \rangle =
(1-A^{-4}) \langle  \sg \rangle
+  A^{2} \langle  )( \rangle =a_1 \langle \sg $\
$\rangle +a_2 \langle )( \rangle $,\
$\langle  T^D \rangle=\langle  \parbox{0.7cm}{
\psfig{figure=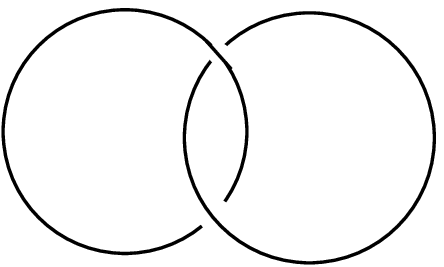,height=0.4cm}}
\rangle=-A^4-A^{-4}$,\ $\langle (T^{(m)})^D \rangle = (-A^4-A^{-4})^m$,
$a_1^m = (1-A^{-4})^m$, and
$a_2^{(m)}= \frac{(-A^4-A^{-4})^m - (1-A^{-4})^m} {d}$.
From this follows that
$\langle  (T^{(m)})^N \rangle=(1-A^{-4})^m(-A^2-A^{-2})+
\frac{(-A^4-A^{-4})^m-(1-A^{-4})^m}{-A^2-A^{-2}}$,\
establishing the first part of Proposition 5.15.\\
To prove Proposition 5.15(ii) we use the formula (iii) of Lemma 5.9
for the product of 2-tangles:\\
$\langle M_{[m[\frac{1}{2}],[s]]} \rangle = a_1^m \langle [s]^N\rangle +
a_2^{(m)}\langle [s]^D\rangle$, and the result of a simple calculation:
$\langle [s]^D\rangle = (-A^3)^s$ and $\langle [s]^N\rangle =
A^{-s}d - A^sd^{-1}(A^{-2s}- (-1)^sA^{2s})$. Thus,
$\langle M_{[m[\frac{1}{2}],[s]]} \rangle = \\
(1-A^{-4})^m\bigl(A^{-s}d - A^sd^{-1}(A^{-2s}- (-1)^sA^{2s})\bigr) +
(-A^3)^s\frac{(-A^4-A^{-4})^m - (1-A^{-4})^m} {d} =$ \\
$(-A^3)^s\bigl(\frac{(-A^4-A^{-4})^m}{d} +
(1-A^{-4})^m\bigl((-A^{-4})^s d -d^{-1}) - (-A^{-2})^s(\frac{A^{-2s} -
(-A^2)^s}{d}\bigr)\bigr)=$\\
$(-A^3)^s\bigl(\frac{(-A^4-A^{-4})^m}{d} +
(1-A^{-4})^m\bigl(-A^{-4})^s(d -d^{-1})
 \bigr)$.

Finally, notice that
for $m\geq 1$ one has  $sw(M_{[m[\frac{1}{2}],[s]]}) = s$ and
therefore $\tilde V_{M_{[m[\frac{1}{2}],[s]]}}(t) =
(-A^3)^{-s}\langle M_{[m[\frac{1}{2}],[s]]} \rangle$, for
$t=A^{-4}$ giving the formula for
$\tilde V_{M_{[m[\frac{1}{2}],[s]]}}(t)$ in Proposition 5.15.
\end{proof}
\begin{proposition}\label{Proposition 5.16}
(i) Pretzel links $M_{[m[\frac{1}{2}],[s]]}$, for $m\geq 3$ are not
$5$-move equivalent to rational links or their connected sums,
and with the additional assumption
that $-2 \leq s \leq 2$, they are pairwise not $5$-move equivalent.\\
(ii) If $m,m'\geq 3$, $m\neq m'$, then
$V(M_{[m[\frac{1}{2}],[s]]}) \neq V(M_{[m'[\frac{1}{2}],[s']]})$, for
any $s,s' \in Z$. Furthermore, $V(M_{[m[\frac{1}{2}],[s]]}) \geq V(T_2)
\approx 1.90211$, and the equality holds only for $m=3, s\equiv 1 \mod 5$;\\
(iii) $V(M_{[m[\frac{1}{2}],[s]]}) = V(M_{[m[\frac{1}{2}],[s']]})$ if
and only if $s \equiv s'\mod 5$ or $m+s+s'\equiv 0 \mod 5$.\\
(iv)
$V_{M_{[m[\frac{1}{2}],[s]]}}(t,5) = V_{M_{[m[\frac{1}{2}],[s']]}}(t,5)$
if and only if $s \equiv s'\mod 5$.
\end{proposition}
\begin{proof} The main tool in our proof is the formula (iii) of
Proposition 5.15:
$$\tilde V_{M_{[m[\frac{1}{2}],[s]]}}(t)=
(-1)^{m-1}\frac{(t+t^{-1})^m}{t^{1/2} + t^{-1/2}} +
(1-t)^m(-t)^s(d - \frac{1}{d}).$$
For $t=e^{\pi i/5}$ (and $t^{1/2}= e^{\pi i/10}$), the first term is a real
number equal approximately to $(-1)^{m-1}\frac{1.618^m}{1.902}$,
independent on $s$ and diverging to infinity.
We can think of this term as the leading
term of the formula\footnote{It is not that unexpected
as $t+t^{-1}=\frac{1+\sqrt 5}{2}\approx 1.61803$
is the golden ratio.}. The second, ``small", term of the formula
 has the absolute value approximately equal to
$0.618^m(1.902 - \frac{1}{1.902}) \approx (1.376)(0.618^m)$, converging
to zero.

 This approximation, with a support of data in Example 5.14,
 suffices to justify (ii) of Proposition 5.16.

A more careful look at the formula for $\tilde V_{M_{[m[\frac{1}{2}],[s]]}}(t)$
allows us to conclude that the leading term not only ``fixes" $m$ but also:\\
(1) If $\tilde V_{M_{[m[\frac{1}{2}],[s]]}}(t) \et
\tilde V_{M_{[m'[\frac{1}{2}],[s']]}}(t)$, $m,m'\geq 3$ then
$\tilde V_{M_{[m[\frac{1}{2}],[s]]}}(t) \stackrel{I_t}{\equiv}
\tilde V_{M_{[m'[\frac{1}{2}],[s']]}}(t)$, and \\
(2)  If $\tilde V_{M_{[m[\frac{1}{2}],[s]]}}(t) \stackrel{I_t}{\equiv}
\tilde V_{M_{[m'[\frac{1}{2}],[s']]}}(t)$,
then $m=m'$ and $s' \equiv s \mod 5$. \\
(3) $V(M_{[m[\frac{1}{2}],[s]]}) = V(M_{[m'[\frac{1}{2}],[s']]})$,
if and only if $m=m'$ and $s' \equiv s \mod 5$ or $m+s+s' \equiv 0 \mod 5$.\\
To see (3) let us rewrite the formula in the form:
$$\tilde V_{M_{[m[\frac{1}{2}],[s]]}}(t)=
(-1)^{m-1}\frac{(t+t^{-1})^m}{t^{1/2} + t^{-1/2}} +
(-t)^{\frac{m}{2} +s}(-t)^{1/2}+(-t^{-1/2})^m(d - \frac{1}{d}).$$
Then it is clear that in order to have $V(M_{[m[\frac{1}{2}],[s]]}) =
V(M_{[m'[\frac{1}{2}],[s']]})$ we need
$\tilde V_{M_{[m'[\frac{1}{2}],[s']]}}(t)$ to be equal to
$\tilde V_{M_{[m[\frac{1}{2}],[s]]}}(t)$ or its conjugate
$\tilde V_{M_{[m[\frac{1}{2}],[s]]}}(t^{-1})$ (all this for
$t^{1/2}= e^{\pi i/10}$). The conjugate condition gives
$m+s+s' \equiv 0 \mod 5$. Finally, as we already noted, the mirror image
of $M_{[m[\frac{1}{2}],[s]]} = M_{[-m[\frac{1}{2}],[-s]]}=
M_{[m[\frac{1}{2}],[-s-m]]}$. Finally, the Jones polynomial of the mirror
image $\bar L$ of $L$ satisfies $\tilde V_{\bar L}(t) = \tilde V_L(t^{-1})$.

This completes the proof of (3) and of Proposition 5.15(iii).
\end{proof}

We will end this section by completing the proof of Theorem 5.1(2).
\begin{proposition}
For $m\geq 3$ a pretzel link $M_{[m[\frac{1}{2}],[s]]}$ is not
$5$-move equivalent to a Montesinos link
$M_{[k[\frac{2}{5}],m[\frac{1}{2}]]}$, $k\geq 1$  and to
a connected sum of $H$'s, $4_1$'s and $T_2$.
\end{proposition}

\begin{proof}
If a link $L$ contains $4_1$ as a summand then $V(L)=0$ and
$L$ cannot be $5$-move equivalent to $M_{[m[\frac{1}{2}],[s]]}$.
If $L$ is a connected sum of
$p$ copies of $H$ and $n$ copies of $T_2$ then
$\tilde V_L(t) = (-t-t^{-1})^p(-t^{1/2}- t^{-1/2})^n \et (1+t^2)^p(1+t)^n$.
Furthermore $L$ can be reduced by $p$ $(2,2)$-moves to the trivial link
of $k+1$ components, thus $col_5(L) = 5^{k+1}$ (also
$F_L(1,2cos(2\pi/5) = (-1)^p(\sqrt 5)^{k})$. On the other hand
$M_{[m[\frac{1}{2}],[s]]}$ can be reduced to the $(2, -2m+s)$ torus
link by $m$ $(2,2)$-moves, thus
\begin{displaymath}
col_5(M_{[m[\frac{1}{2}],[s]]}) =
\left\{
\begin{array}{lll}
Z_5^2 & for   & -2m+s \equiv 0\ \mod 5\\
Z^5 & for & -2m+s \not\equiv 0\ \mod 5
\end{array}
\right .
\end{displaymath}
Therefore, in the connected sum we can have only none  or one copy
of $T_2$. Consider this two cases independently:\\
($k=0$) Then we would have $-2m+s \equiv 0\ \mod 5$ and for
$\sqrt t= e^{\pi/10}$:
$$\frac{(t+t^{-1})^m}{t^{1/2} + t^{-1/2}} +
(-1)^{m-1}((-t)^{1/2}+(-t^{-1/2}))^m(d - \frac{1}{d}) =
(t+t^{-1})^p.$$
It would imply that $m+1 \leq p \leq m+1$, the contradiction.\\
(k=1) Then, as in the case of $k=0$ we are forced
 have $-2m+s \equiv 0\ \mod 5$
and for $\sqrt t= e^{\pi/10}$:\\
$$\frac{(t+t^{-1})^m}{t^{1/2} + t^{-1/2}} +
(-1)^{m-1}((-t)^{1/2}+(-t^{-1/2}))^m(d - \frac{1}{d}) =
(t+t^{-1})^p(t^{1/2} + t^{-1/2}),$$ or equivalently
$$(t+t^{-1})^m + (-1)^{m-1}((-t)^{1/2}+(-t^{-1/2}))^m(1-d^2) =
(t+t^{-1})^p(t+t^{-1} +1),$$
which is impossible.

To complete the proof of Proposition 5.16 we should distinguish
a pretzel link $M_{[m[\frac{1}{2}],[s]]}$, $m\geq 3$
from a Montesinos link
$M_{[k[\frac{2}{5}],n[\frac{1}{2}]]}$, $k\geq 1, k+n\geq 3$.
The consideration is
similar to the previous one. First we use Fox 5-coloring to see
that for Montesinos links it is $5^k$ or $5^{k+1}$ while for a pretzel
knot $5$ or $5^2$ thus $k=1$ or $2$ and $n\geq 1$. Then we use
our formulas for the Jones polynomial
and comparing their values for $\sqrt t = e^{\pi/10}$ we see that the
right side has a real representative (when considered up to $\pm t^{i/2}$,
therefore the left side should have a real representative, which forces us
to have $-2m+s \equiv 0\ \mod 5$. With this, we would have the equality:
$$\frac{(t+t^{-1})^m}{t^{1/2} + t^{-1/2}} +
(-1)^{m-1}((-t)^{1/2}+(-t^{-1/2}))^m(1-d^2) =
|(1-t)^{k+n-2} (1+t)^{k-1}|.$$
We quickly find it impossible for $k=1$ or $k=2$. Namely, for  $k\leq2$
we have
$|(1-t)^{k+n-2} (1+t)^{k-1}|\leq |(1-t)(1+t)| \approx 1.175$, which is
smaller then the possible values for the left hand sight (which, as
we already computed is $\geq |1+t|\approx 1.902$). The proof of
Proposition 5.16 is complete.
\end{proof}

\ \\

\section{Density of values of $V(L) = |V_L(e^{\pi i/5})|$}\label{Section 6}

We recover here V.~F.~R.~Jones observation (\cite{Jon}, Corollary 14.7) that
the values of $V(L)=|V_L(e^{\pi i/5})|$ are dense in $[0,\infty)$.
Furthermore, Lemma 6.1 has played an important role in our proof
of Corollary 5.13.

\begin{lemma}\label{Lemma 6.1}\
 Let $k_1,k_2,{k_1}',{k_2}' \geq 0$, then \\
(i) $(1+t)^{k_1}(1-t)^{k_2} \et 1 $ if and only if $k_1=k_2=0$\\
~~(ii) $(1+t)^{k_1}(1-t)^{k_2} \et (1+t)^{{k_1}'}(1-t)^{{k_2}'}$
if and only if $k_1={k_1}'$ and $k_2={k_2}'$
\end{lemma}
~~\\
\begin{proof}
 First we note that (ii) follows from (i). Namely, without lost of
 generality we can assume $k_1 \geq {k_1}',$
 if $k_2 < {k_2}'$ then from the fact that $1+t$ and
 $1-t$ are not zero divisors in $Z[t]/({\frac{1+t^5}{1+t}})$, we have
 $(1+t)^{{k_1}-{k_1}'}
 \et (1-t)^{{k_2}-{k_2}'}$.\\
For $t=e^{\pi i /5}$ we would have \\
$1 \leq \mid 1+ e^{\pi i/5} {\mid}^{k_1-{k_1}'} =
\mid 1- e^{\pi i /5} \mid < 1$ the contradiction.\\
In proving (i) first assume that $k_1 \geq 1$ and
consider the equation in (i)
modulo ideal $(\frac{t^5+1}{t+1},t+1)=(t+1,5)$,
then we have $0 \equiv \pm t^{i}$ mod $(t+1,5)$, the contradiction. \\
Therefore we have $k_1=0$ and $(1-t)^{k_2} \et 1$.
Then $\mid 1-e^{\frac{\pi i}{5}} {\mid}^{k_2}=1 $
which holds only for $k_2=0$. \\
The proof of Lemma 6.1 is completed.
\end{proof}

\begin{corollary}\label{Corollary 6.2}
The values of $\mid V_{K}(e^{\frac{\pi i}{5}}) \mid $
for $K$ a connected sum of any number of copies of knots
$5_1$ and $8_{17}$ form a dense subset of $(0, \infty)$.
\end{corollary}
\begin{proof}
$V_{5_1}(t) \et 1+t$ and $V_{8_{17}}(t) \et 1-t$
therefore by Lemma 6.2 the values \\
$\mid V_{{k_1}{5_1} \# {k_2}{8_{17}} }(e^{{\pi i}/5})
\mid =\mid 1+e^{\pi i/5} {\mid}^{k_1}
\mid 1-e^{{\pi i}/5} {\mid}^{k_2}$ are all different\\
 and  never equal to $1$.\\
Because $1 < {\mid} 1+e^{\pi i/5} {\mid} \approx 1.90211$
and   $1 > {\mid} 1-e^{\pi i/5} {\mid} \approx 0.618034$
therefore the values ${\mid} 1+e^{{\pi i}/5} {\mid}^{k_1}
 {\mid} 1-e^{{\pi i}/5} {\mid}^{k_2}$ taken over all positive
 $k_1$ and $k_2$ are dense in $(0, \infty)$. \\
\end{proof}

\begin{corollary}\label{Corollary 6.3}
The values of
$\mid V_{M_{[k[\frac{2}{5}],m[\frac{1}{2}] ]}}(e^{\pi i/5}) \mid$
for $k,m >0$, form a dense subset of $(0, \infty)$.
It is the case because in the formula of Theorem 4.6,
$1 < {\mid} 1+e^{2\pi i/5} {\mid} \approx 1.17557$
and   $1 > {\mid} 1-e^{\pi i/5} {\mid} \approx 0.618034$

\end{corollary}

We can interpret Corollary 6.3 as suggesting that classification
of links up to $5$-moves is as difficult as classification of
links in general. However, the goal of this paper was to show
that for some classes classification is to some degree possible.
Motivated by the case of $(2,2)$-moves we had in mind the class
of algebraic links. The classification of rational and pretzel links and
the partial classification of Montesinos links
is the first step in this direction.

\section{Tables of links up to $9$-crossings}\label{Section 7}
In the following table we list all prime links up to 9 crossings and
some of their $5$-move invariants. In our notation, $r$ before the
name of a link denotes rational link, $p$ denotes non-rational
pretzel link and $m$ denotes a Montesinos link which is neither rational
nor pretzel link. ${}^*$ before the name of a link denotes 
a link which is not $5$-move equivalent to its mirror image. The letter
$a$ after the name of the knot denotes
an amphicheiral knot. Links in the same ``box" are $5$-move equivalent. 
If the representative of a box (in the first column) is in the Bold 
face then the links in the box are not $5$-move equivalent to links 
in any other box.

~~\\
{\footnotesize{
\begin{center}
\begin{tabular}{|l|l|l|l|l|l|} \hline
 Rep. &
${\bf{F}}$& $V_{L,{\overline{L}}}(t,5)$ & $L$  & ${\overline{L}}$ & V \\    \hline
  {\bf{T}}$_1$ &
  $1$  &
$1$ &$rT_1 $,$r6_1$,$r6_2$,$r7_2$,$r7_6$ && 1\\
 & &
& $r7_7$,$r8_4$,$r8_{12}a$,$r8_{13}$ & & \\
  & &
& $r8_{14}$,$r9_{1}$,$r9_{3}$,$r9_{4}$& & \\
 & &
& $r9_{7}$,$r9_{8}$,$r9_{9}$, $r9_{15}$& & \\
& &
& $r9_{17}$,$r9_{18}$,$r9_{19}$, $r9_{20}$& & \\& &
& $r9_{27}$,$r4_1^2$,$r6_1^2$, $r7_1^2$& & \\
 & &
& $r7_3^2$,$r8_2^2$,$r8_4^2$,$r8_5^2$ & & \\
 & &
& $r8_8^2$,$r9_4^2$,$r9_6^2$,$r9_7^2$& & \\
& &
& $r9_8^2$,$r9_{10}^2$, $r9_{11}^2$, $p9_{19}^2$& & \\
& &
&  $m9_{20}^2$,$p9_{22}^2$, $p9_{51}^2$, $p8_7^3$ &  & \\
& &
&  $p9_5^3$, $m9_6^3$ &  & \\
\hline
{\bf{T}}$_2$  &
$\sqrt{5}$ &$1+t$ & $r5_1 \stackrel{5}\sim T_2$ &  & 1.90211\\   &
$$&&  $r7_4$,$r8_8$,$r9_{23}$,$r9_{31}$ & & \\
 & $$&& $r8_6^2$,$r8_7^2$,$r9_9^2$,$r9_{12}^2$ & & \\
 &
$$& &  $p9_{50}^2$, $p8_1^3$ &  & \\
\hline
{\bf{4}}$_1$ &
 $-\sqrt{5}$&  $0$  &  $r4_1a$,$r8_9a$,$r9_2$,$r9_{12}$  & & 0 \\
 & &
    &  $r6_2^2$, $r9_1^2$,$r9_3^2$&  & \\ \hline
 {\bf{H}} &$-1 $ &
$1+t^2$ & $H$, $r3_1$, $r5_2$, $r6_3a$ & &1.61803 \\
& $ $ &
  & $r7_1$,$r7_3$,$r7_5$, $r8_1$, $r8_2$& &
  \\
     & $$ &
  & $r8_3a$, $r8_6$, $r8_7$,$r8_{11}$& & \\
     & $$ &
 &$r9_5$,$r9_6$,$r9_{10}$,$r9_{11}$, $r9_{13}$ & & \\
  & $ $ &
 & $r9_{14}$, $r9_{21}$,$r9_{26}$, $r2_1^2$  &  & \\
   & $ $ &
 & $r5_1^2$,$r6_3^2$,$r7_2^2$,$r8_1^2$,$r8_3^2$&  & \\
 &&&$r9_2^2$,$r9_5^2$, $p9_{21}^2$, $9_{37}^2$ &  &\\
 &&& $p9_{49}^2$,$p9_{52}^2$, $p8_2^3$,$p8_8^3$&  &\\
   \hline
{\bf{8}}$_{21}$  & $-\sqrt{5}$ & $1+t$&  $p8_{21} \stackrel{(2,2)}\sim T_2$ &  & 1.90211\\
 &
$$& & $p9_{24}$, $p9_{37}$, $p7_5^2$, $m9_{15}^2$ &  & \\
 &
$$ & &  $p9_{27}^2$, $9_{34}^2$, $m9_{48}^2$, $p9_{53}^2$ &  & \\
 &
$$& &  $p7_1^3$, $m9_{14}^3$ &  & \\
\hline
{\bf{9}}$_{17}^3$ &$ 1$ &
$1+2t+t^2-t^3$ & ${}^*p9_{17}^3$ &  ${}^*p9_3^3$,  ${}^*p8_2^4$& 3.44298\\
\hline
{\bf{9}}$_3^3$ &$ 1 $&
$1-t-2t^2-t^3$ & ${}^*p9_3^3$, ${}^*p8_2^4$ &  ${}^*p9_{17}^3$ &3.44298 \\
\hline
${\bf{8}}_{14}^2$ &$1$ &
$2+t$ &   ${}^*8_{14}^2$, ${}^*9_{33}^2$, ${}^*9_{10}^3$ & ${}^*9_{56}^2$, ${}^*9_{18}^3$ & 2.86986\\ \hline
${\bf{9}}_{56}^2$  &
$1$ &$1+2t$ & ${}^*9_{56}^2$, ${}^*9_{18}^3$ &  ${}^*8_{14}^2$, ${}^*9_{33}^2$, ${}^*9_{10}^3$ &2.86986 \\
\hline
 \end{tabular}
\newpage
\begin{tabular}{|l|l|l|l|l|l|} \hline
 Rep. &
${\bf{F}}$& $V_{L,{\overline{L}}}(t,5)$ & $L$  & ${\overline{L}}$ & V \\    \hline
{\bf{9}}$_{49}$ &
$-5$&$1+t$&  $9_{49} \stackrel{(2,2)}{\nsim} T_2$  &  & 1.90211\\
\hline
{\bf{8}}$_{18}$ &$\sqrt{5}$ &
$2-t+2t^2$ &  $8_{18}a$ &  &$\sqrt{5}$ \\   \hline
{\bf{9}}$_{20}^3$ &$-1$  &
$1-2t$ &  ${}^*9_{20}^3$ & & 1.32813\\  \hline
$\overline{{\bf{9}}_{20}^3}$ & $-1$& $1-2t+t^2$  & & ${}^*9_{20}^3$ &1.32813 \\   \hline
${\bf{9}}_{15}^3$ & $\sqrt{5}$ &
$2+2t$ &  $p9_{15}^3$, $p8_3^4$ & &3.80423 \\
\hline
${\bf{9}}_{13}^2$ &$1$ &
$1+t+t^2$ &$p9_{13}^2$,$p9_{14}^2$ & & \\
 & &
&   $p9_{43}^2$,  $p9_{44}^2$, $H \# H$& &2.61803 \\
\hline
  {\bf{9}}$_{12}^3$  & 1&
$3-t+2t^2$ &  ${}^*9_{12}^3$ & & 3.1013 \\   \hline
${\overline{{\bf{9}}_{12}^3}}$ & $1$ &
$2+t-t^2-t^3$ & &  ${}^*9_{12}^3$ & 3.1013\\    \hline
{\bf{9}}$_{40}$ & 5&
$1-2t+2t^2$ &  $9_{40}$ & & 0.726543\\
  \hline
${\bf{9}}_{31}^2$
&$-\sqrt{5}$ &
$1+t+t^2+t^3$ & $9_{31}^2$ & & 3.07768\\
 \hline
{\bf{8}}$_{10}$ & $-1$&
$1+t+t^2-t^3$ & ${}^*p8_{10}$, ${}^*p7_8^2$, ${}^*p8_{11}^2$  & ${}^*p8_{15}$,${}^*p8_{19}$, ${}^*p9_{35}$   &
2.49721\\
& $$ &
&  ${}^*m9_{45}^2$,${}^*p9_{54}^2$, ${}^*p6_1^3$  & ${}^*m9_{48}$,${}^*p8_{12}^2$, ${}^*m9_{18}^2$  & \\
 & $$ &
&& ${}^*m9_{47}^2$, ${}^*p8_3^3$,${}^*8_5^3$  & \\
&&&& ${}^*m9_2^3$, ${}^*m9_{13}^3$&\\
\hline
{\bf{8}}$_{15}$&$ -1$&
$1-t-t^2-t^3$ &${}^*p8_{15}$,${}^*p8_{19}$, ${}^*p9_{35}$ &  ${}^*p8_{10}$, ${}^*p7_8^2$, ${}^*p8_{11}^2$ &
2.49721\\
   &  &
 & ${}^*m9_{48}$,${}^*p8_{12}^2$, ${}^*m9_{18}^2$ &  ${}^*m9_{45}^2$,${}^*p9_{54}^2$, ${}^*p6_1^3$& \\
&&& ${}^*m9_{47}^2$, ${}^*p8_3^3$,${}^*8_5^3$&&\\
&&& ${}^*m9_2^3$, ${}^*m9_{13}^3$&&\\
\hline
{\bf{9}}$_9^3$ & $-1$&
$1+t+t^2-2t^3$ & ${}^*9_9^3$ & & 2.76008\\  \hline
${\overline{{\bf{9}}_9^3}}$ &
$ -1$ &
$2-t+3t^2-2t^3$ & &  ${}^*9_9^3$ & 2.76008\\    \hline
{\bf{9}}$_{21}^3$
& $-\sqrt{5}$ &
$2-3t+2t^2-3t^3$ & ${}^*9_{21}^3$& & 2.93565\\    \hline
${\overline{{\bf{9}}_{21}^3}}$ &
 $-\sqrt{5}$ &
$3-2t+3t^2-2t^3$ &&  ${}^*9_{21}^3$& 2.93565\\   \hline
 {\bf{9}}$_{55}^2$
 & $1 $ &
$2-2t-2t^3$ & $9_{55}^2$,  $6_2^3$,$8_6^3$ & & 3.23607 \\ \hline
{\bf{8}}$_5$ & $1$ &
$2-t+2t^2-t^3$ & ${}^*p8_5$, ${}^*m9_{28}$, ${}^*p9_{46}$ & ${}^*p8_{20}$, ${}^*p9_{16}$, ${}^*p7_4^2$ & 2.14896
\\
 & $ $ &
& ${}^*p7_7^2$,${}^*m9_{17}^2$, ${}^*m9_{24}^2$&  ${}^*7_6^2$,${}^*m9_{16}^2$, ${}^*p9_{23}^2$  & \\
      & $ $ &
& ${}^*p9_{28}^2$, ${}^*9_{32}^2$,${}^*9_{35}^2$ &${}^*9_{29}^2$, ${}^*p6_3^3$ & \\
  & $ $ &
&${}^*9_{39}^2$, ${}^*m9_{46}^2$, ${}^*9_{59}^2$   & & \\
&&& ${}^*8_9^3$,  ${}^*m9_1^3$, ${}^*m9_7^3$ && \\  \hline
{\bf{8}}$_{20}$ &$ 1$
$ $&
$1-2t+t^2-2t^3$ &  ${}^*p8_{20}$, ${}^*p9_{16}$, ${}^*p7_4^2$ & ${}^*p8_5$, ${}^*m9_{28}$, ${}^*p9_{46}$  &
2.14896 \\
  & $ $ &
&  ${}^*7_6^2$,${}^*m9_{16}^2$, ${}^*p9_{23}^2$  &  ${}^*p7_7^2$,${}^*m9_{17}^2$, ${}^*m9_{24}^2$ & \\
  & $ $ &
& ${}^*9_{29}^2$, ${}^*p6_3^3$ & ${}^*p9_{28}^2$, ${}^*9_{32}^2$,${}^*9_{35}^2$  & \\
& $ $ &
&&${}^*9_{39}^2$, ${}^*m9_{46}^2$, ${}^*9_{59}^2$  & \\
&&&& ${}^*8_9^3$,  ${}^*m9_1^3$, ${}^*m9_7^3$&\\
\hline
\end{tabular}
\newpage
\begin{tabular}{|l|l|l|l|l|l|} \hline
Rep. &
${\bf{F}}$&$V_{L,{\overline{L}}}(t,5)$ & $L$  & ${\overline{L}}$  & V \\   \hline
$9_{38}$ & $-1$&
$2$ &  ${}^*9_{38}$, ${}^*9_{30}^2$,  ${}^*9_8^3$ &  &2 \\ \hline
$9_{47}$& $-1$&
$2$ &   ${}^*9_{47}$ &   &2 \\ \hline
${\overline{9_{38}}}$ & $-1$&
$2$ &  &  ${}^*9_{38}$, ${}^*9_{30}^2$,  ${}^*9_8^3$ &2 \\ \hline
${\overline{9_{47}}}$ & $-1$&
$2$ &   & ${}^*9_{47}$  &2 \\ \hline
$8_{17}$   &
$-1$ &
$1-t$ &  $8_{17}a$, $9_{58}^2$, $8_4^3$ & & 0.618034 \\  \hline
$9_{22}$ & $-1$ &  $1-t$
&  $m9_{22}$,  $m9_{25}$,$m9_{30}$ & & 0.618034\\  &  &
& $m9_{36}$, $m9_{42}$,  $m9_{43}$ & & \\
 &  &
& $m9_{44}$, $m9_{45}$, $m8_9^2$ & & \\
& $$ &
& $m8_{10}^2$,$m8_{15}^2$, $m8_{16}^2$  & & \\
&&&$m9_{25}^2$, $m9_{26}^2$,$9_{36}^2$ & & \\
\hline
$9_{41}^2$&$1$ &
$1+t-t^2$ & ${}^*9_{41}^2$& ${}^*9_{42}^2$& 1.54336 \\ \hline
$9_{60}^2$&$1$ &
$1+t-t^2$ & ${}^*9_{60}^2$, ${}^*9_{19}^3$ & ${}^*9_{29}$& 1.54336 \\ \hline
${\overline{9_{32}}}$&$1$ &
$1+t-t^2$ &  & ${}^*9_{32}$, ${}^*9_{33}$& 1.54336 \\ \hline
${\overline{9_{41}}}$&$1$ &
$1+t-t^2$ &  & ${}^*9_{41}$& 1.54336 \\ \hline

$9_{29}$ & $1$  &
$1-t-t^2$ &  ${}^*9_{29}$ & ${}^*9_{60}^2$, ${}^*9_{19}^3$&1.54336 \\ \hline
$9_{32}$ & $1$  &
$1-t-t^2$ &  ${}^*9_{32}$,  ${}^*9_{33}$ & &1.54336 \\ \hline
${\overline{9_{41}^2}}$ & $1$  &
$1-t-t^2$ & ${}^*9_{42}^2$ & ${}^*9_{41}^2$&1.54336 \\ \hline
$9_{41}$ & $1$  &
$1-t-t^2$ &  ${}^*9_{41}$ & &1.54336 \\ \hline
$9_{34}$       &$1$ &
$1-2t+t^2$ &  ${}^*9_{34}$ &   & 0.381966\\  \hline
${\overline{9_{34}}}$       &$1$ &
$1-2t+t^2$ &   &  ${}^*9_{34}$ & 0.381966\\  \hline
$8_{16}$
& $\sqrt{5}$&
$1-t+2t^2-t^3$ & ${}^*8_{16}$, ${}^*9_{57}^2$ &  ${}^*8_{10}^3$& 1.17557\\ \hline
$8_{10}^3$
& $\sqrt{5}$&
$1-t+2t^2-t^3$ & ${}^*8_{10}^3$ &  ${}^*8_{16}$,  ${}^*9_{57}^2$& 1.17557\\ \hline
$8_{13}^2$
& $\sqrt{5}$&
$1-t+2t^2-t^3$ & $8_{13}^2$, $9_{38}^2$, $9_{11}^3$ & & 1.17557\\ \hline
$9_{39}$ & $-\sqrt{5}$ &$1-t+2t^2-t^3$
& ${}^*9_{39}$ &   & 1.17557\\ \hline
${\overline{9_{39}}}$ & $-\sqrt{5}$ &$1-t+2t^2-t^3$
&  &  ${}^*9_{39}$ & 1.17557\\ \hline
$9_{40}^2$ & $5 $ &$1-t+2t^2-t^3$& ${}^*9_{40}^2$ & ${}^*9_{61}^2$ & 1.17557\\ \hline
$9_{61}^2$ & $5 $ &$1-t+2t^2-t^3$& ${}^*9_{61}^2$ &  ${}^*9_{40}^2$ & 1.17557\\ \hline
$9_4^3$ &
 $ -1$ &
$2+2t^2+t^3$&  $m9_4^3$, $p9_1^4$ &  & 3.67044\\  \hline
${\overline{9_{16}^3}}$ &
 $ -1$ &
$2+2t^2+t^3$&  &  ${}^*m9_{16}^3$, ${}^*p8_1^4$ & 3.67044\\  \hline
$9_{16}^3$& $-1 $ &
$1+2t+2t^3$ & ${}^*m9_{16}^3$, ${}^*p8_1^4$ &  &3.67044 \\  \hline
${\overline{9_4^3}}$& $-1 $ &
$1+2t+2t^3$ &  &  ${}^*m9_4^3$, ${}^*p9_1^4$ &3.67044 \\  \hline\end{tabular}
\end{center}
}}

\section{Acknowledgment}
M. Ishiwata was supported by the 21st COE program
``Constitution of wide-angle mathematical basis focused on knots".
J. Przytycki was supported by JSPS grants when visiting Japan
in December 2004 -January 2005 and June 2007.

\end{document}